\documentclass[10pt, a4paper]{article}
\author{Frank Schuhmacher\footnote{Supported by:
Doktorandenstipendium des Deutschen Akademischen Austauschdienstes im Rahmen
des gemeinsamen Hochschulsonderprogramms III des Bundes und der
L\"ander}}

\title{Deformation of $L_\infty$-Algebras}

\usepackage{mathenv}
\usepackage{a4wide}
\usepackage[intlimits]{amsmath}
\usepackage{amssymb}
\usepackage{amsthm} 
\usepackage{amscd}
\usepackage[all]{xypic}
\usepackage{makeidx}

\newcounter{punkt}

\theoremstyle{definition}
\newtheorem{defi}{Definition}[section]
\newtheorem{bem}[defi]{Remark}
\newtheorem{prop}[defi]{Proposition}

\newtheorem{beisp}[defi]{Example}

\theoremstyle{theorem}
\newtheorem{satz}[defi]{Theorem}

\newtheorem{lemma}[defi]{Lemma}
\newtheorem{kor}[defi]{Corollary}

\newcommand{\nach}{\longrightarrow}

\newcommand{\isom}{\cong}

\newcommand{\NN}{ \mathbb{N} }
\newcommand{\RR}{ \mathbb{R} }
\newcommand{\ZZ}{\mathbb{Z}}

\newcommand{\ppp}{\cdot\ldots\cdot}

\newcommand{\im}{\operatorname{im}}

\newcommand{\I}{\operatorname{I}}

\renewcommand{\S}{\operatorname{S}}

\newcommand{\bew}{\bf Proof: \rm}
\renewcommand{\qed}{\hfill $ \square$ \\}

\newcommand{\Hom}{\operatorname{Hom}}

\newcommand{\id}{\operatorname{Id}}

\newcommand{\pr}{\operatorname{pr}}

\newcommand{\Kern}{\operatorname{Kern} }

\newcommand{\Kokern}{\operatorname{Cokern} }

\newcommand{\lz}{\hfill\newline}

\newcommand{\ot}{\otimes}
\newcommand{\lie}{[\cdot,\cdot]}
\newcommand{\Coder}{\operatorname{Coder}}
\newcommand{\odos}{\odot\ldots\odot}
\newcommand{\ots}{\otimes\ldots\otimes}
\newcommand{\sig}{\sigma}
\newcommand{\Sh}{\operatorname{Sh}}
\newcommand{\Ot}{\operatorname{Ot}}
\newcommand{\ddd}{,\ldots,}
\newcommand{\kkk}{+\ldots+}
\newcommand{\Wedge}{\bigwedge}
\newcommand{\pnu}{\downarrow}
\newcommand{\pno}{\uparrow}

\newcommand{\dgmanf}{\operatorname{DG-Manf}}
\newcommand{\Def}{\operatorname{Def}}

\begin{document}

\maketitle
\begin{abstract}
In this paper, deformations of $L_\infty$-algebras 
are defined in such a way that the bases of deformations are
$L_\infty$-algebras, as well. A
universal and a semiuniversal deformation is constructed for
$L_\infty$-algebras, whose cotangent complex admits a splitting.
The paper also contains an
explicit construction of a minimal $L_\infty$-structure
on the homology $H$ of a differential graded Lie algebra $L$
and of an $L_\infty$-quasi-isomorphism between $H$ and $L$.  
\end{abstract}

\tableofcontents
\pagebreak
\section*{Introduction}

$L_\infty$-algebras (see Section~\ref{general})
play a crucial role in deformation 
theory. They are natural generalizations of differential
graded Lie algebras (DGLs).
Deformation problems can always
be described by DGLs (see \cite{Kont}, for instance). 
The importance of $L_\infty$-algebras in deformation theory
comes from the fact that
two different deformation problems are equivalent, if
the corresponding DGLs are equivalent as $L_\infty$-algebras.
This was one ingredient of Kontsevich's \cite{Kont} prove that
deformation quantization works on each Poisson manifold.
$L_\infty$-algebras also build a bridge from algebra to
geometry. A simple shift of degrees makes a formal DG-manifold
out of an $L_\infty$-algebra (see Section~\ref{general}).
This observation is also due to Kontsevich.
If a deformation problem is governed by a DGL $L$, then
the (formal) local muduli space, if it exists, is an analytic
subspace of the formal DG-manifold corresponding to
$L$.\\

In the other way, to each DGL $L$, one can define an abstract 
deformation functor $\Def_L$. In the classical theory
$\Def_L$ is a set-valued functor on the category of
(Artinian, local) algebras. Recent studies in
mirror symmetry (\cite{KontS}, \cite{Merk2}) 
led to an extension of this functor first to
graded and then differential graded Artinian algebras.
The aim of this extension is to produce smooth (in a sense) formal
moduli spaces with tangent space, isomorphic to
the whole cohomology of $L$. 
But sometimes it is not evident (or not even possible) to
give an algebraic or geometric meaning to the objects obtained
by the extended deformation functor. 
(Sometimes this is possible. 
For the classical deformations of
associative algebras, the extended deformation functor produces
the more general $A_\infty$-algebras.)\\

The deformation theory of $L_\infty$-algebras
(or in geometric
terms, of formal DG-manifolds) presented in this paper
is in fact an extended deformation theory of (formal) singularities. 
Instead of working with deformation functors, we
present a completely geometric (extended) deformation
theory of formal DG manifolds. The bases of deformations are
formal DG manifolds as well. The theory is developed 
analogous to ``embedded deformations'' of singularities.
The deformations of a given formal DG-manifold $M=(M,Q^M)$
are governed by the DGL $L$ of formal vectorfields on $M$
(see Section~\ref{defolinf}). The degree 1 shift of $L$
is again a formal DG manifold, denoted by $U$. There are
two nice observations. The first is that the 
going over $M\mapsto U$ doesn't change the category. (This one is
trivial.) The second (Theorem~\ref{unidef}) is that $U$ is the base
of a universal deformation of $M$.   
For the construction of a semiuniversal deformation of
$M$, we have to construct an $L_\infty$-structure on
the homology $H$ of $(L,d)$, such that $H$ and $L$ are
equivalent as $L_\infty$-algebras. 
$H$ with such an $L_\infty$-structure is called a minimal
model for $L$.\\

Hence, the essence of this paper
is the following general recipe for the construction of (formal)
analytic moduli spaces: 
Take a minimal representative in the class of $L_\infty$-algebras
modulo $L_\infty$-equivalence of
the DGL controlling the deformation problem.
This recipe had been discovered before (see \cite{KontL},
\cite{Merk1}) and was rediscovered independently by the author.\\

The contents of this paper:
In Section~\ref{general} we remind the definitions 
of $L_\infty$-algebras and of their correspondence with differential graded
coalgebras. We state the conditions for a sequence of maps,
to define an $L_\infty$-morphism. We will prove those conditions in
the Appendix, since they are hard to find in the literature.
Then we remind Kontsevich's geometric point of view (=formal DG manifolds)
of $L_\infty$-algebras.
In Section~\ref{defolinf}, we define deformations of formal DG
manifolds
with formal DG bases and morphisms of those. 
Our definition generalizes the one of Fialowski and Penkava \cite{Penk}.
We show that for an arbitrary formal DG manifold $M$, the differential
graded Lie algebra $\Coder(S(M),S(M))$ (which we call tangent
complex of $M$) is a base of a universal deformation of $M$.
In Section~\ref{trees}, we give an ad-hoc combinatorial introduction
to binary trees. In a sense, binary trees contain the algebraic
structure of $L_\infty$-algebras (see \cite{Voro}). In
Section \ref{leq}, they are used to define an $L_\infty$-structure $\mu_\ast$
on the homology $H$ of a differential graded Lie algebra $L=(L,d,\lie)$ 
(admitting a splitting).
Furthermore, again with the help of binary trees, in Section~\ref{leq},
we construct explicitly an
$L_\infty$-quasi-isomorphism between $(H,\mu_\ast)$ and $(L,d,\lie)$.
Similar constructions in the $A_\infty$-context are due to
Gugeheim/Stasheff \cite{GuSt}, Merkulov
\cite{MerkK} and Kontsevich/Soibelman \cite{KontS}.
In Section~\ref{hodgecomp}, we prove that $(L,d,\lie)$ is 
as $L_\infty$-algebra
isomorphic to the direct sum of the $L_\infty$-algebras
$(H,\mu_\ast)$ and $(F,d)$,
where $F$ is the complement of $H$ in $L$.
As consequence, we can show that for each formal DG manifold
$M$ such that
$L=\Coder(S(M),S(M))$ splits, the shift $V$ of
$(H,\mu_\ast)$ is the base of
a semi-universal deformation of $M$.\\

\textbf{Acknowledgements:} I want to express my
gratitude to Siegmund Kosarew for inspiring this
work and for many valuable suggestions and discussions.

\section{$L_\infty$-Algebras and Coalgebras}\label{general} 

In this paper we shall always work over a ground ring $k$
of characteristic zero.

\subsection{Graded Symmetric and Exterior Algebras}

For a graded module $W$, the graded symmetric algebra $S(W)$ is
defined as the tensor algebra $T(W)=\oplus_{n\geq 0}W^{\ot n}$
modulo the relations $w_1\ot w_2-(-1)^{w_1 w_2}w_2\ot w_1=0$.
We denote the graded symmetrical product by $\odot$.
The algebra $T(W)$ (resp.$S(W)$) is bigraded. The graduation on $T(W)$
(resp. $S(W)$), defined by 
$g(w_1\ot...\ot w_n)=g(w_1)+...+g(w_n)$ 
(resp. $g(w_1\odot...\odot w_n)=g(w_1)+...g(w_n)$), where
$g$ is the graduation of $W$, will be called \textbf{linear graduation}.
The one defined by $g(w_1\ot...\ot w_n)=n$ 
(resp.$g(w_1\odot...\odot w_n)=n$), 
will be called \textbf{polynomial graduation}.
Set $S_+(W):=\oplus_{n\geq 1}W^{\odot n}$.
On $S_+(W)$, there is a natural $k$-linear comultiplication
$\Delta^+:S_+(W)\nach S_+(W)\ot S_+(W)$, given by
$$w_1\odot...\odot w_n\mapsto\sum_{j=1}^{n-1}\sum_{\sig\in\Sh(j,n)} 
\epsilon(\sigma,w_1,...,w_n)w_{\sig(1)}\odos w_{\sig(j)}\ot
w_{\sig(j+1)}\odos w_{\sig(n)}.$$
Here $\epsilon(\sig):=\epsilon(\sigma,w_1,...,w_n)$ is defined such that
$w_{\sig(1)}\odos w_{\sig(n)}=\epsilon(\sig)w_1\odos w_n$.
Note that we have $\Kern\Delta^+=W$.
\newcommand{\uw}{\underline{w}}
On $S(W)$, there is a $k$-linear comultiplication $\Delta$, defined by
$\Delta(1):=1\ot 1$ and $\Delta(\uw):=\uw\ot 1+\Delta'(\uw)+1\ot \uw$,
for $\uw\in\S_+(W)$. Note that $\Delta$ is injective.

For a graded module $L$, the graded exterior algebra $\Wedge^+ L$ 
without unit is
defined as the tensor algebra $T_+(L)=\oplus_{n\geq 1}L^{\ot n}$
modulo the relations $a_1\ot a_2+(-1)^{a_1 a_2}a_2\ot a_1=0$.
We denote the graded exterior product by $\wedge$.
$L[1]$ denotes the graded module with $L[1]^i=L^{i+1}$
and $\pnu$ the canonical map $L\nach L[1]$ of degree $-1$.
Set $\pno:=\pnu^{-1}$.
Remark that, for each $n\geq 1$, there is an isomorphism
\begin{align*}
\pnu^n:{\Wedge}^n L&\nach L[1]^{\odot n}\\
a_1\wedge...\wedge a_n & \mapsto (-1)^{(n-1)\cdot a_1\kkk 1\cdot
  a_{n-1}} \pnu a_1\odos\pnu a_n
\end{align*}

Its inverse map is given by $(-1)^{\frac{n(n-1)}{2}}\pno^n$.
As we shall always do, we just have applied the Koszul sign
convention, i.e. for homogeneous graded morphisms $f,g$
of graded modules, we set $(f\ot g)(a\ot b):=(-1)^{ga}f(a)\ot g(b)$. 
In the exponent, $a$ always means the degree
of an homogeneous element (or morphism) $a$ and $ab$ means
the product of degrees, not the degree of the product.\\

For $\sig\in\Sigma_n$ and $a_1\ddd a_n\in L$, we define the
sign $\chi(\sig):=\chi(\sig,a_n\ddd a_n)$ in such a way that
$$a_{\sig(1)}\wedge...\wedge a_{\sig(n)}=\chi(\sig)a_1\wedge...\wedge a_n.$$
We have the following correlation between $\chi$ and $\epsilon$:
\begin{bem}
For $a_1,...,a_n\in L$, we have
$$\chi(\sig,a_1\ddd a_n)=
(-1)^{(n-1)(a_1+a_{\sig(1)})\kkk 1\cdot(a_{n-1}+a_{\sig(n-1)})}
\epsilon(\sig,\pnu a_1\ddd\pnu a_n).$$
\end{bem}
For a graded module $V$,
we define two different actions of the symmetric group
$\Sigma_n$ on $V^{\ot n}$:
The first one is given by
\begin{align*}
\Sigma_n\times V^{\ot n}&\nach V^{\ot n}\\
(\sig,v_1\ot...\ot v_n)&\mapsto\epsilon(\sig,v_1,...,v_n)
v_{\sig(1)}\ot...\ot v_{\sig(n)}
\end{align*}
Here, the application of a $\sig$
commutes with the canonical projection
$V^{\ot n}\nach V^{\odot n}$.
The second one is given by
\begin{align*}
\Sigma_n\times V^{\ot n}&\nach V^{\ot n}\\
(\sig,v_1\ot...\ot v_n)&\mapsto\chi(\sig,v_1,...,v_n)
v_{\sig(1)}\ot...\ot v_{\sig(n)}
\end{align*}
Here, the application of a $\sig$
commutes with the canonical projection
$V^{\ot n}\nach \wedge^nV$.
When we work with symmetric powers, we use the first
action; when we work with exterior powers, we use the
second one. Since the context shall always be clear, 
we don't distinguish both actions by different notations.
We will use the anti-symmetrisation maps:
$$\alpha_n:=\sum_{\sig\in\Sigma_n}\sig.:V^{\ot n}\nach V^{\ot n}.$$
When $\sig.$ denotes the first action, $\alpha_n$
can be seen as map $V^{\odot n}\nach V^{\ot n}$; when $\sig.$
denotes the second action, $\alpha_n$
can be seen as map $\wedge^n V\nach V^{\ot n}$.
Furthermore, for both cases, we will use the maps
$$\alpha_{k,n}:=\sum_{\sig\in\Sh(k,n)}\sig.:
V^{\ot n}\nach V^{\ot n}.$$

For the natural projection $\pi:W^{\ot n}\nach W^{\odot n}$
(resp. $V^{\ot n}\nach \Wedge^nV$), we have
$$\pi\circ\alpha=n!\id.$$

\subsection{Free Differential Graded Coalgebras}

Let $(C_1,\Delta_1)$ and $(C_2,\Delta_2)$ be coalgebras.
Remember that a module homomorphism 
$F:C_1\nach C_2$ is a coalgebra map,
iff the diagram

\begin{equation}\label{qudiag}
\xymatrix{
C_1\ar[r]^F\ar[d]^{\Delta_1} & C_2\ar[d]^{\Delta_2}\\
C_1\ot C_1\ar[r]^{F\ot F} & C_2\ot C_2
}\end{equation}

commutes. Each coalgebra morphism
$F:(S(W),\Delta)\nach (S(W'),\Delta)$ satisfies $F(1)=1$. 
The restriction $F\mapsto F|_{S_+(W)}$
gives is a 1:1-correspondence
between coalgebra morphisms
$(S(W),\Delta)\nach (S(W'),\Delta)$ and coalgebra morphisms
$F:(S_+(W),\Delta^+)\nach (S_+(W'),\Delta^+)$.

The next proposition gives a one-to-one correspondence
between coalgebra maps $ F:S(W)\nach S(W')$
and sequences of linear maps 
$F_n:S_n(W)\nach W'$, $n\geq 1$.
We fix the following notations:
\begin{align*}
 \hat{F}_n:=&  F|_{W^{\odot n}}: W^{\odot n}\nach  S(W')\\
F_{k,l}:=& \pr_{{W'}^{\odot l}}\circ\hat{F}_k: W^{\odot k}
\nach {W'}^{\odot l}\\
F_n:=&F_{n,1}: W^{\odot n}\nach W'
\end{align*}

Sometimes, we shall consider the maps $F_n$ as antisymmetric
maps $W^{\ot n}\nach W$ instead of maps $W^{\odot n}\nach W$. 
For each multi-index $I=(i_1,...,i_k)\in\NN^k$, we set
$I!:=i_1!\ppp i_k!$ and $|I|:=i_1\kkk i_k$ and
$$F_I:=\frac{1}{I!k!}(F_{i_1}\odot...\odot F_{i_k})\circ\alpha_n.$$

Here, by $F_{i_1}\odot...\odot F_{i_k}$, we mean the composition
of $F_{i_1}\ot...\ot F_{i_k}$ and the natural projection
${W'}^{\ot k}\nach {W'}^{\odot k}$.

\begin{prop}\label{Fcoal}
For $n\geq 1$, we have that
\begin{equation}
 \hat{F}_n=\sum_{k=1}^n\sum_{I\in\NN^k}^{|I|=n}
F_I.
\end{equation}
\end{prop}
The proof can be found in the appendix.\\

A coalgebra homomorphism $F:S(W)\nach S(W')$ is called \textbf{strict},
if $F_n=0$ for each $n\geq 2$.\\ 

For a coalgebra $(C,\Delta)$,
remember that a module homomorphism $Q:C\nach C$ is a coderivation,
iff the diagram
\newcommand{\hQ}{\hat{Q}}

\begin{equation}\label{qdelta}
\xymatrix{
C\ar[rr]^{Q}\ar[d]^\Delta && C\ar[d]^{\Delta}\\
C\ot C\ar[rr]^{Q\ot 1+1\ot Q}& & C\ot C
}\end{equation}

commutes.
By the next proposition, there is a one-to-one correspondence
between coderivations $Q:S(W)\nach S(W)$ of degree $+1$ and 
sequences of linear maps $Q_n:S_n(W)\nach W$ of degree $+1$.
We fix the following notations:

\begin{align*}
\hQ_n:=&Q|_{W^{\odot n}}: W^{\odot n}\nach S(W)\\
Q_{k,l}:=&\pr_{W^{\odot l}}\circ \hQ_k: W^{\odot k}\nach W^{\odot l}\\
Q_n:=&Q_{n,1}: W^{\odot n}\nach W
\end{align*}

\begin{prop}\label{Qcoal} Let $Q$ be a coderivation of degree $+1$
on the coalgebra $(S(W),\Delta)$. Then, $Q(1)=Q_0(1)\in W$ and
for $n\geq 1$ and $w_1,...,w_n\in W$, we have
\begin{equation}
\hQ_n(w_1,...,w_n)=\sum_{l=0}^n\sum_{\sig\in\Sh(l,n)}
\epsilon(\sig)
Q_l(w_{\sig(1)}\ddd w_{\sig(l)})\odot
w_{\sig(l+1)}\odos w_{\sig(n)}, 
\end{equation}
where the $l=0$ term must be interpreted as
$Q_0(1)\odot w_1\odos w_n$.
\end{prop}

The proof can be found in the appendix.
Remark that there is a 1:1-correspondence between
coderivations of degree $+1$ on $(S_+(W),\Delta^+)$
and coderivations of $Q$ degree $+1$ on $(S(W),\Delta)$ with
$Q(1)=0$.

\begin{kor}
Let $Q$ be a coderivation of degree $+1$ on the coalgebra $S(W)$, 
$Q'$ a coderivation of degree $+1$ on the coalgebra $S(W)$ and
$ F:=S(W)\nach S(W')$ a morphism of coalgebras. Then, for $n\geq 1$ and
$1\leq l\leq n+1$, $1\leq k\leq n$,  we have
\begin{equation*}
Q_{n,l}=
(Q_{n-l+1}\ot 1\ot...\ot 1)\circ\alpha_{n-l+1,n}.
\end{equation*}
and
\begin{equation*}
F_{n,k}=\sum_{i_1\kkk i_k=n}F_I.
\end{equation*}
$F$ respects the coderivations $Q$ and $Q'$ iff 
$F(Q(1))=Q'(1)$ and for each $n\geq 1$
we have
\begin{equation}\label{coalgmor}
\sum_{k=1}^n\sum_{I\in\NN^k\atop |I|=n}
Q'_k\circ F_I=
\sum_{k+l=n+1}
F_l\circ(Q_{k}\ot 1\ots 1)\circ\alpha_{k,n}
\end{equation}
On the right hand - side, the sum is over all $l\geq1$ and $k\geq 0$.
$(Q_0\ot 1\ots 1)(w_1\ots w_n)$ must be interpreted as
$Q_0(1)\ot w_1\ots w_n$.
\end{kor}

\subsection{$L_\infty$-Algebras}

Remember that a module $L$ with a sequence of maps 
$\mu_n:\Wedge^nL\nach L$ of degree $2-n$, for $n\geq 0$, is called an
\textbf{$L_\infty$-algebra} iff the coderivation $Q$
(of degree $+1$) on $S(W)$, defined by the maps 
$$Q_n:=(-1)^{n(n-1)/2}\pnu\circ\mu_n\circ\pno^n:W^{\odot n}\nach W$$ 
is a codifferential, i.e. $Q^2=0$.
\begin{bem}
The condition $Q^2=0$ just means that for each $n$ and $w_1\ddd
w_n\in W$ the term
$$(Q^2)_n(w_1\ddd w_n)=
\sum_{k+l=n+1}\sum_{\sig\in\Sh(k,n)}\epsilon(\sig)
Q_l(Q_{k}(w_{\sig(1)}\ddd w_{\sig(k)}), w_{\sig(k+1)}\ddd 
w_{\sig(n)})
$$
disappears.
This conditions can easily be translated in the following ones
\begin{align}
\sum_{k+l=n+1}\sum_{\sig\in\Sh(k,n)}(-1)^{k(l-1)}\chi(\sig)
\mu_l(\mu_{k}(a_{\sig(1)}\ddd a_{\sig(k)}), a_{\sig(k+1)}\ddd 
a_{\sig(n)})=0
\end{align}
for each $n\geq 0$ and $a_1\ddd a_n\in L$.
\end{bem}

In the literature $\mu_0$ is mostly assumed to be trivial.
If this is the case, $(L,\mu_1)$ is a DG module.

\begin{defi}\label{mini}
An $L_\infty$-algebra $(L,\mu_n)_{n\geq 1}$ is called \textbf{minimal}, if
$\mu_1=0$. It is called \textbf{linear}, if
$\mu_i=0$ for $i\geq 2$.
\end{defi}

Now let $(L,\mu_\ast)$ and $(L',\mu'_\ast)$ be $L_\infty$-algebras.
Set $W:=L[1]$, $W':=L'[1]$ and denote the induced codifferentials on
$S(W)$ and $S(W')$ by $Q$ and $Q'$.
A sequence of maps
$f_n:\Wedge^nL\nach L';\;n\geq 0$ of degree $1-n$ is called
\textbf{$L_\infty$-morphism} iff the maps $F_n:=W^{\odot n}\nach W'$ 
induced by $f_n$ (explicitly: 
$F_n=(-1)^{n(n-1)/2}\pnu\circ f_n\circ\pno^n$) define a map
$ F:S(W)\nach S(W')$ of differential graded coalgebras. 
Rewrite
condition (\ref{coalgmor}) into terms of $f_n$ and $\mu_n$:

$$\begin{array}{c}
\mu'_1\circ f_n-\sum_{i+j=n}
\frac{(-1)^i}{2}\mu'_2(f_i,f_j)\circ\alpha_{i,j}
+\sum_{k=3}^n\sum_{I\in\NN^k\atop |I|=n}
(-1)^{k(k-1)/2+i_1(k-1)+...+i_{k-1}\cdot 1}\mu'_k\circ f_I\\
\quad\\
=\sum_{k+l=n+1}(-1)^{k(l-1)}f_l\circ(\mu_k\ot 1\ot...\ot 1)\circ\alpha_{k,n}
\end{array}$$

For the case
where $L'$ is a differential graded Lie-algebra, i.e.
$\mu'_k=0$ for $k=0$ and $k\geq 3$, set
$d:=\mu'_1$ and $[\;,\;]:=\mu'_2$. Then we get the following
conditions for the maps $f_n$ to define an $L_\infty$-morphism
(see Definition 5.2 of \cite{Lama}):

$$\begin{array}{l}
df_n(a_1\ddd a_n)-\\
\quad\\
-\sum_{i+j=n}\sum_\sigma\chi(\sig)
(-1)^{i+(j-1)(a_{\sig(1)}\kkk a_{\sig(i)})}[f_i(a_{\sig(1)}\ddd
a_{\sig(i)}),f_j(a_{\sig(i+1)}\ddd a_{\sig(n)})]\\
\quad\\
\quad\quad=\sum_{k+l=n+1}\sum_{\sig\in\Sh(k,n)}
(-1)^{k(l-1)}\chi(\sig)f_l(\mu_k(a_{\sig(1)}\ddd
a_{\sig(k)}),a_{\sig(k+1)}\ddd a_{\sig(n)}),
\end{array}$$

where $a_1,...,a_n\in L$ and
the second sum goes over all $\sig$ in $\Sh(i,n)$
such that $\sig(1)<\sig(i+1)$.

\begin{defi}
A morphism $f:(L,\mu_n)_{n\geq 1}\nach (L',\mu'_n)_{n\geq 1}$ of 
$L_\infty$-algebras is called 
\textbf{$L_\infty$-quasi-isomorphism},
if $f_1$ is a quasi-isomorphism of differential graded modules.
\end{defi}

\subsection{$L_\infty$-Algebras and Formal DG Manifolds} \label{super}

In this subsection, we explain briefly the geometric point of view
of $L_\infty$-algebras, as proposed by Kontsevich \cite{Kont}.
First, recall the definition of \textbf{pointed modules} (see
Section II.6 of \cite{BiKo}): A pointed module is a pair $(M,\ast)$ of
a module $M$ and an element $\ast\in M$. We restrict ourselves to
the case where $\ast$ is just the zero element.
For modules $M$ and $N$,
a \textbf{homogeneous polynomial
of degree} $p$ on $M$ with values in $N$
is a mapping $\tilde{f}:M\nach N$ of the form
$f\circ{\Delta}^{(p)}$, where $f$ is a $p$-multilinear form
$M\times...\times M\nach N$ and ${\Delta}^{(p)}$ is the
diagonal $m\mapsto(m,...,m)$.
The polarization formula (Lemma II.6.2 of \cite{BiKo}) says
that $f\mapsto\tilde{f}$ is a
1:1-correspondence between symmetrical $p$-multilinear forms
$M\times...\times M\nach N$ and
homogeneous polynomials of degree $p$
on $M$ with values in $N$.
For pointed modules $M=(M,0)$ and $N=(N,0)$, a formal map
$f:M\nach N$ is a formal sum $f=\sum_{p\geq 1}\tilde{f}_p$,
where
$\tilde{f}_p$ is a homogeneous polynomial of degree $p$.
Pointed modules, together with formals maps form a category.
By the polarization formula and 
Proposition~\ref{Fcoal}, there is
a 1:1-correspondence between
formal maps $f:M\nach N$ and morphisms $S(M)\nach S(N)$ of
(non-graded free) coalgebras.\\

For the definition of formal supermanifolds, we replace modules
by $\ZZ$-graded modules and symmetric multilinear forms by
graded symmetric multilinear forms.

\begin{defi}
A \textbf{formal  supermanifold} is a pair $M=(M,0)$ of
a $\ZZ$-graded module $M$ and its zero element.
A \textbf{formal map} $f:M\nach N$ of degree $j$ of formal supermanifolds
is a sequence $(f_p)_{p\geq 1}$, where $f_p$ is a graded symmetric
multilinear form $M\times...\times M\nach N$ of linear degree $j$.
The \textbf{composition} $f\circ g$ of formal maps $g:L\nach M$ and
$f:M\nach N$ is defined as the sequence $(g_p)_{p\geq 1}$
with 
$$g_p=\sum_{k=1}^p\sum_{I\in\NN^k\atop |I|=p}f_k\circ g_I.$$ 
A \textbf{morphism} of formal supermanifolds is a formal map
of degree zero.
\end{defi}
It is clear by this definition that the category of formal
supermanifolds is equivalent to the category of
free, graded coalgebras with coalgebra maps of degree zero.

\begin{defi}
A \textbf{vectorfields} of degree $j$ on a formal supermanifold
$M$ is a coderivation of degree $j$ on $S(M)$.
\end{defi}

By Proposition~\ref{Qcoal}, 
a vectorfields on $M$ can be interpreted as
formal map $M\nach M$. The $M$ on the right hand-side of the arrow should be
considered as tangent space of $M$. 
The graded commutator defines the structure
of a graded Lie algebra on $\Coder(S(M),S(M))$. Therefore, there
is a bracket $\lie$ of vectorfields.\\

Let $(M,Q^M)$ and $(N,Q^N)$ be formal supermanifolds with vectorfields.
A formal map $f:M\nach N$ is called \textbf{$Q$-equivariant}, if
the induced map $S(M)\nach S(N)$ of coalgebras commutes with
$Q^M$ and $Q^N$.
Remark that in the case where $M$ and $N$ are non-graded free
$k$-modules of finite dimension, this definitions coincide with
the classical definitions and the $Q$-equivariance just means that
$$Q^N\circ f=Df\circ Q^M.$$

\begin{defi}
A \textbf{formal DG manifold} is a pair $(M,Q^M)$ of a formal
supermanifold $M$ and a vectorfield $Q^M$ of degree $1$ such that
$[Q^M,Q^M]=0$. Morphisms of formal DG manifolds are $Q$-equivariant
maps of formal supermanifolds (sometimes we call them
$L_\infty$-morphisms).
Denote the category of formal DG manifolds by
$\dgmanf$. 
\end{defi}

By the previous subsection, the lifting $L\mapsto L[1]$ gives a
1:1-correspondence between $L_\infty$-algebras and
formal DG manifolds
and the functor $M\mapsto S(M)$ gives a 1:1-correspondence
between formal DG manifolds and free differential graded
coalgebras.

We use the following superscripts to denote full
subcategories of $\dgmanf$:\\
L (``local''): the subcategory of all $(M,Q^M)$ in $\dgmanf$
such that $Q^M_0=0$;\\
M (``minimal''): the subcategory of all $(M,Q^M)$ in $\dgmanf^L$
such that $Q^M_1=0$;\\
G (``g-finite''): the subcategory of all $(M,Q^M)$ in $\dgmanf^L$ such that 
$H(M,Q^M_1)$ is g-finite;\\
C (``convergent''): the subcategory of all $(M,Q^M)$ in 
$\dgmanf^GM$ such that the mapping $M_0\nach M_1$ induced by $Q^M$ 
converges.\\

\section{Deformation of $L_\infty$-Algebras}\label{defolinf}

Fialowski and Penkava \cite{Penk} have defined a
deformation theory of $L_\infty$-algebras such that the base of
a deformation is an algebras with augmentation. 
The new approach here is 
to take $L_\infty$-algebras also as bases of deformations.
Since the geometric language is more elegant, we will talk about
formal DG manifolds instead of $L_\infty$-algebras. Thus, the
objects that we deform are DG structures, i.e.
degree 1 vectorfields $Q$ with $Q^2=0$ on
formal supermanifolds.\\

In our setting, not every ``fiber'' of a deformation of 
a DG structure on $M$ gives DG structure 
on $M$ but in general
only a degree 1 vectorfield. 
But it is easy to find those points of the basis $B$ of a deformation of
$M$ for which the
associated deformation of $Q^M$ is again a DG structure.
They just correspond to the zero locus of the vectorfield $Q^B$.\\

A very nice fact for this deformation theory is that
we get a universal deformation for free: The deformations of
a DG manifold $M$ are governed by the differential
graded Lie algebra of vectorfields on $M$, 
i.e. the DGL $L$ of coderivations on $S(M)$ with
graded commutator as bracket $\lie$ and differential
$d=[\cdot,Q^M]$.
In contrast to Fialowski/ Penkava, we use the linear grading
on $L$ (see Section~\ref{general}). 
Set $U:=L[1]$ and denote the vectorfield corresponding to the DGL structure
of $L$ by $Q^U$. We will see that $(U,Q^U)$ is the base of a universal 
deformation  of $M$. 

\subsection{Definitions}

\begin{defi}
Let $(M,Q^M)$ in $\dgmanf$ and $(B,Q^B)$ in
$\dgmanf^L$ be formal DG manifolds. A \textbf{deformation}
of $M$ with base $B$, or more exactly a deformation of the DG structure
$Q^M$, is a degree 1 vectorfield $Q$ on $B\times M$ with $Q_0=0$, such that
\begin{enumerate}
\item
$Q|_{\{0\}\times M}=0$.   
\item
$\tilde{Q}:=Q^M+Q^B+Q$ is a DG structure on $B\times M$.
\item
The projection $B\times M\nach B$ is a homomorphism of
formal DG manifolds.
\end{enumerate}
\end{defi}
We denote deformations of $(M,Q^M)$ as triples $(B,Q^B,Q)$.
Remark that condition (i) is equivalent to the condition that
the inclusion $M\nach B\times M$ is a
morphism of formal DG manifolds.  Condition (iii) is equivalent
to the condition
$$\im(Q)\subseteq\{0\}\times M.$$
A deformation is trivial, if the projection
$B \times M\nach M$ respects the DG structures
$\tilde{Q}$ and $Q^M$.
\begin{defi}
A \textbf{morphism} of deformations $(B,Q^B,Q)$ and 
$(B',Q^{B'},Q')$ of $(M,Q^M)$ is a pair $(F,f)$,
where $F$ is a morphism of formal DG manifolds
$(B\times M,\tilde{Q}:=Q^B+Q^M+Q)$ and 
$(B'\times M,\tilde{Q}':=Q^B+Q^M+Q')$ and $f$ is a morphism of
formal DG manifolds $(B,Q^B)$ and $(B',Q^{B'})$ 
such that the diagram 
$$\xymatrix{
B\times M\ar[r]\ar[d] & B'\times M\ar[d]\\
B\ar[r] & B'}$$
is cartesian and the diagram
$$\xymatrix{
M\ar[d]\ar[dr] & \\
B\times M\ar[r] & B'\times M
}$$
commutes.
\end{defi}

\begin{defi}
Two deformations are called \textbf{equivalent}, if there exist
homomorphisms in both senses.
\end{defi}

\begin{prop}\textbf{(Base change)}
Suppose that $(B',Q^{B'},Q')$ is a deformation of $(M,Q^M)$ and
$f:B\nach B'$ a homomorphism of formal DG manifolds with
$B=(B,Q^B)$ in $\dgmanf^L$.
Then, via
$$Q_n(b_1,...,b_r,m_1,...,m_s):=
\sum_{t=1}^r\sum_{\I\in\NN^t\atop |I|=r}
Q'_{s+t}(f_I(b_1,...,b_r),m_1,...,m_s)$$
for $r\geq 1$ with $r+s=n$ and $b_i\in B'$, $m_j\in M$, we
can define a deformation $(B,Q^{B},Q)$ of $(M,Q^M)$ 
and $(f\times\id,f)$ is a morphism of deformations.
\end{prop}
\bew
We have to show that $(Q^2)_n(b_1,...,b_r,m_1,...,m_s)=0$, for
$b_1,...,b_r\in B$ and $m_1,...,m_s\in M$. First let
$s\geq 1$. Then
\begin{multline*}
(Q^2)_n(b_1,...,b_r,m_1,...,m_s)=\\
\sum_{k=r}^n\sum_{\sig\in\Sh(k-r,s)}\epsilon(\sig)
Q^M_{n-k+1}(Q_k(b_1,...,b_r,m_{\sig(1)},...,m_{\sig(k-r)}),
m_{\sig(k-r+1)},...,m_{\sig(s)})+\\
\sum_{k=1}^r\sum_{\sig\in\Sh(k,r)}\epsilon(\sig)
Q_{n+1-k}(Q_k^B(b_{\sig(1)},...,b_{\sig(k)}),
b_{\sig(k+1)},...,b_{\sig(r)},m_1,...,m_s)+\\
\sum_{k=0}^s\sum_{\sig\in\Sh(k,s)}\epsilon(\sig)
Q_{n+1-k}(b_1,...,b_r,Q^M_k(m_{\sig(1)},...,m_{\sig(k)}),m_{\sig(k+1)},
...,m_{\sig(s)})+\\
\sum_{k+l=n+1}\sum_{p=1}^{r-1}\sum_{\sig\in\Sh(p,r)}
\sum_{\tau\in\Sh(k+p-r,s)}
\epsilon(\sig)\epsilon(\tau)
Q_l(b_{\sig(1)},...,b_{\sig(p)},Q_k(b_{\sig(p+1)},...,b_{\sig(r)},\\
m_{\tau(1)},...,m_{\tau(k+p-r)}),m_{\sig(k+p-r+1)},...,m_{\tau(s)}).
\end{multline*}

Using the definition of $Q$ and the assumption that
$f$ is a DG morphism, 
after changing the order of summation, this sum takes the form
\begin{multline*}
\sum_{t=1}^r\sum_{I\in\NN^t\atop |I|=r}\sum_{p=0}^s\sum_{\sig\in\Sh(p,s)}
\epsilon(\sig)Q^M_{s-p+1}(Q'_{p+t}(f_I(b_1,...,b_r),m_{\sig(1)},
...,m_{\sig(p)}),m_{\sig(p+1)},...,m_{\sig(s)})+\\
\sum_{p=1}^r\sum_{t=1}^p\sum_{I',I'',u}\sum_{\sig\in\Sh(p,r)}\epsilon(\sig)
Q'_{s+u+1}(Q_t^B(f_{I'}(b_{\sig(1)},...,b_{\sig(p)})),
f_{I''}(b_{\sig(p+1)},...,b_{\sig(r)}),m_1,...,m_s)+\\
\sum_{t=1}^r\sum_{I\in\NN^t\atop |I|=r}\sum_{p=0}^s\sum_{\sig\in\Sh(p,s)}
\epsilon(\sig)Q'_{s-p+t+1}(f_I(b_1,...,b_r),Q^M_p(m_{\sig(1)},...,m_{\sig(p)}),
m_{\sig(p+1)},...,m_{\sig(s)})+\\
\sum_{p=1}^{r-1}\sum_{q=0}^s
\sum_{t=1}^p\sum_{I',I'',u}\sum_{\sig\in\Sh(p,r)}\sum_{\tau\in\Sh(q,s)}
\epsilon(\sig)\epsilon(\tau)
Q'_{t+s-q+1}(f_{I'}(b_{\sig(1)},...,b_{\sig(p)}),\\[2mm]
Q'_u(f_{I''}(b_{\sig(p+1)},...,b_{\sig(r)}),m_{\tau(1)},...,m_{\tau(q)}),
m_{\tau(q+1)},...,m_{\tau(s)}),
\end{multline*}
where in the second and forth term, the sum is taken over all
$I'\in\NN^t$ such that $|I'|=p$, over all $u=1,...,r-p$
and all $I''\in\NN^u$ such that $|I''|=r-p$.
But this sum equals
$$\sum_{t=1}^r\sum_{I\in\NN^t\atop |I|=r}
(\tilde{Q})^2_{s+t}(f_I(b_1,...,b_r),m_1,...,m_s),$$
which is zero. The case $s=0$ goes in the same manner.\\

Now, $f\times\id$ is a map of formal DG manifolds
and for a diagram
$$\xymatrix{
A\ar@/_/[ddr]_g\ar@/^/[drr]^h\ar@{.>}[dr] &&\\
& B\times M\ar[r]^{f\times\id}\ar[d] & B'\times M\ar[d]\\
& B\ar[r]^f & B'
}$$

one can show that $j:=g\times(\pr_M\circ h)$ is a
DG morphism completing the diagram commutatively.
Hence, the quadratic diagram is cartesian and the pair $(f\times\id,f)$
is a morphism of deformations.
\qed

\begin{kor}
If $(F,f)$ is a morphism $(B,Q^B,Q)\nach(B',Q^{B'},Q')$
of deformations and $f$ an isomorphism, then $(F,f)$ is also an
isomorphism.
\end{kor}
\bew
The deformation $(B,Q^B,Q)$ is natural isomorphic to the deformation,
obtained by base change. For the latter one, the statement is clear.
\qed

\subsection{A Universal Deformation}

\begin{defi}
A deformation $(U,Q^U,Q)$ of $(M,Q^M)$ is called \textbf{universal}, if for
each deformation $(B,Q^B,Q')$ there exists 
a morphism $(F,f):(B,Q^B,Q')\nach(U,Q^U,Q)$, where
$f$ is uniquely defined. 
A deformation $(V,Q^V,Q)$ is called \textbf{semi-universal}, if for
each deformation $(B,Q^B,Q')$ there exists a homomorphism
$(B,Q^B,Q')\nach(V,Q^V,Q)$ of deformations and if
$(V,Q^V)$ is minimal (in the sense of Definition~\ref{mini}).
\end{defi}

Let $L$ be the differential graded Lie algebra $\Coder(S(M),S(M))$
with bracket
$$[s,t]=s\circ t-(-1)^{st}t\circ s,$$
for homogeneous $s,t$
and differential $d(s):=(-1)^s[s,Q^M]$.

\begin{defi}
The complex $(L,d)$ will be called \textbf{tangent complex} of $M$.
\end{defi}

Set $U:=L[1]$ and denote the vectorfield corresponding to the DGL structure
on $L$ by $Q^U$. There is a canonical construction of a deformation $Q$
of $M$ with base $U$:\\ 

Define multilinear maps
$q_n:U\ot M^{\ot n-1}\nach M$ of degree $+1$ by
$$u\ot m_1\ot...\ot m_{n-1}\mapsto (\uparrow u)(m_1\odos m_{n-1})$$ 
and denote the symmetrisation of the map  
$(U\times M)^{\ot n}\nach U\times M$, induced by $\frac{1}{n!}q_n$ by $Q_n$.
Hence, we get a vectorfield $Q$ of degree $+1$ on $U\times M$, such that
$Q|_{\{0\}\times M}=0$. 
 
\newcommand{\tQ}{\tilde{Q}}
\begin{prop}
$\tilde{Q}:=Q^M+Q^U+Q$ is an DG vectorfield on $U\times M$ and
the projection $U\times M\nach U$ respects the DG structures
$\tilde{Q}$ and $Q^U$.
\end{prop}

\bew
Remember that we have
$$(\tilde{Q}^2)_n=\sum_{k+l=n+1}\tQ_l\circ(\tilde{Q}_k\ot1\ot...\ot1)
\circ\alpha_{k,n}.$$
Since $Q^M$ and $Q^U$ are $L_\infty$-structures, we have
$(\tQ^2)_n(a_1,...,a_n)=0$ if all $a_i$ belong to $M$ or
if all $a_i$ belong to $U$. Hence, it is enough to show that
$(\tQ^2)_n$ is zero on products of the form
\begin{enumerate}
\item[(a)]
$w\odot m_2\odos m_n$, for $n\geq 2$, $w\in U$ and 
$m_2,...,m_n\in M$,
\item[(b)]
$w_1\odot w_2\odot m_3\odos m_n$, for $n\geq 3$, $w_1,w_2\in U$ and 
$m_3,...,m_n\in M$.
\end{enumerate}
We use the abbreviations $m'_i:=m_{i+1}$ and $m''_i:=m_{i+2}$.
For products of the form (a), we have
\begin{spliteqn*}
(\tQ^2)_n(w,m_2,...,m_n)=\\[2mm]
\sum_{k+l=n+1\atop k>1}\sum_{\sig\in\Sh(k,n)\atop \sig(1)=1}
\epsilon(\sig,w,m_2,...,m_n) Q_l^M(Q_k(w,m_{\sig(2)},...,m_{\sig(k)}),
m_{\sig(k+1)},...,m_{\sig(n)})\\
+\sum_{k+l=n+1\atop l>1}\sum_{\sig\in\Sh(k,l)\atop \sig(1)=k+1}
\epsilon(\sig,w,m_2,...,m_n) Q_l(Q_k^M(m_{\sig(1)},...,m_{\sig(k)}),
w,m_{\sig(k+2)},...,m_{\sig(n)})\\
+Q_n(Q^U_1(w),m_2,...m_n) =\\
\quad\\
-\sum_{k'+l=n}\sum_{\sig'\in\Sh(k',n-1)}
\epsilon(\sig',m'_1,...,m'_{n-1}) Q_l^M
(\uparrow w(m'_{\sig'(1)},...,m'_{\sig'(k')}),
m'_{\sig(k'+1)},...,m'_{\sig(n-1)})\\
-\sum_{k+l'=n}\sum_{\sig'\in\Sh(k,n-1)}(-1)^w
\epsilon(\sig',m'_1,...,m'_{n-1}) 
(\uparrow w)(Q_k^M(m'_{\sig'(1)},...,m'_{\sig'(k)}),
m'_{\sig(k+1)},...,m'_{\sig(n-1)})\\
+(Q^M\circ\uparrow w+(-1)^w\uparrow w\circ Q^M)(m'_1,...,m'_{n-1})\\=0.
\end{spliteqn*}

For products of the form (b), we have 
\begin{spliteqn*}
(\tQ^2)_n(w_1,w_2,m_3,...,m_n)\quad=\quad
Q_{n-1}(Q^U_2(w_1,w_2),m_3,...,m_n)+\\[3mm]
\sum_{k+l=n+1\atop k\geq2}\sum_{\sig\in\Sh(k,n)\atop \sig(1)=1,\sig(k+1)=2}
\epsilon(\sig,w_1,w_2,m_3,...,m_n) Q_l(Q_k(w_1,m_{\sig(2)},...,m_{\sig(k)}),
w_2,m_{\sig(k+2)},...,m_{\sig(n)})+\\
\sum_{k+l=n+1\atop k\geq 2}\sum_{\sig\in\Sh(k,n)\atop \sig(1)=2}
\epsilon(\sig,w_1,w_2,m_3,...,m_n) Q_l(Q_k(w_2,m_{\sig(2)},...,m_{\sig(k)}),
w_1,m_{\sig(k+2)},...,m_{\sig(n)})
\end{spliteqn*}
For a $\sig\in\Sh(k,n)$ such that $\sig(1)=1$ (resp. $\sig(1)=2$)
and $\sig(k+1)=2$ (resp. $\sig(k+1)=1$), 
define $\sig'\in\Sh(k-1,n-2)$ by
$\sig'(1)=\sig(2)-2,...,\sig'(k-1):=\sig(k)-2$ and
$\sig'(k):=\sig(k+2)-2,...,\sig'(n-2)=\sig(n)-2$.
Then, we have
$$\epsilon(\sig,w_1,w_2,m_3,...,m_n)=
(-1)^{w_2(m_{\sig(2)}+...+m_{\sig(k)})}\epsilon(\sig',m''_1,...,m''_{n-2})$$
(resp.
$$\epsilon(\sig,w_1,w_2,m_3,...,m_n)=
(-1)^{w_1(m_{\sig(2)}+...+m_{\sig(k)})}
\epsilon(\sig',m''_1,...,m''_{n-2}).\quad)$$
Hence, the above sum takes the form
\begin{spliteqn*}
-(-1)^{w_1}(\uparrow w_1\circ\uparrow w_2)(m_3,...,m_n)
-(-1)^{w_2+w_1w_2}(\uparrow w_2\circ\uparrow w_1)(m_3,...,m_n)\\[2mm]
+\sum_{k'+l'=n-1}\sum_{\sig'\in\Sh(k',n-2)}
(-1)^{w_2(w_1+1)}\epsilon(\sig')
(\uparrow w_2)((\uparrow w_1)(m''_{\sig'(1)},...,m''_{\sig'(k')}),
m''_{\sig(k'+1)},...,m''_{\sig'(n-2)})\\
+\sum_{k'+l'=n-1}\sum_{\sig'\in\Sh(k',n-2)}(-1)^{w_1}
\epsilon(\sig') 
(\uparrow w_1)((\uparrow w_2)(m''_{\sig'(1)},...,m''_{\sig'(k')}),
m''_{\sig(k'+1)},...,m''_{\sig(n-2)})\\
=0,
\end{spliteqn*}
where $\epsilon(\sig')$ stands for $\epsilon(\sig',m''_1,...,m''_{n-2})$.
\qed

Hence, $(U,Q^U,Q)$ is a deformation of $(M,Q^M)$. 

\begin{satz}\label{unidef}
$(U,Q^U,Q)$ is a universal deformation of $Q^M$.
More precisely, the mapping $Q'\mapsto f$, where
$$(\uparrow f_n(b_1\odos b_n))_k(m_1,...,m_k):= 
Q'_{n+k}(b_1,...,b_n,m_1,...,m_k)$$
defines a 1:1-correspondence between deformations of $M$
with base $(B,Q^B)$ and morphisms $B\nach U$
of formal DG manifolds.
\end{satz}

\bew
We have to show that 
$(Q^M+Q^B+Q')^2=0$, iff
the family $(f_n)_n$ defines a map 
$f:S(B)\nach S(U)$ of differential graded
coalgebras, i.e. iff for each $n$, and
$b_1,...,b_n\in B$, the equation
\begin{multline}\label{resc}
Q^U_1(f_n(b_1,...,b_n))+
\frac{1}{2}\sum_{i+j=n}\sum_{\sig\in\Sh(i,n)}\epsilon(\sig,b)
Q^U_2(f_i(b_{\sig(1)},...,b_{\sig(i)})\ot f_j(b_{\sig(i+1)},...,b_{\sig(n)}))
=\\
\sum_{k+l=n+1}\sum_{\sig\in\Sh(k,n)}\epsilon(\sig,b)
f_l(Q^B_k(b_{\sig(1)},...,b_{\sig(k)}),
b_{\sig(k+1)},...,b_{\sig(n)})
\end{multline}
holds.
In equation (\ref{resc}),
we apply both sides on terms $m_1\odos m_r\in M^{\odot r}$ and
use the definition of $f$. Then, the condition on $f$ is
equivalent to the condition that the following term is zero:
\begin{spliteqn*}
\sum_{k+l=r+1}\sum_{\tau\in\Sh(k,r)}\epsilon(\tau,m)
Q_l^M(Q'_{n+k}(b_1,...b_n,m_{\tau(1)},...,m_{\tau(k)}),
m_{\tau(k+1)},...,m_{\tau(r)})+\\
(-1)^{b_1+...+b_n}\sum_{k+l=r+1}\sum_{\tau\in\Sh(k,r)}\epsilon(\tau,m)
Q'_{n+l}(b_1,...b_n,Q^M_k(m_{\tau(1)},...,m_{\tau(k)}),
m_{\tau(k+1)},...,m_{\tau(r)})\\
+\sum_{i+j=n\atop i,j\geq 1}\sum_{\sig\in\Sh(i,n)}\epsilon(\sig,b)
(-1)^{b_{\sig(1)}+...+b_{\sig(i)}}
\sum_{k=0}^n\sum_{\tau\in\Sh(k,r)}\epsilon(\tau,m)
Q'_{i+l}(b_{\sig(1)},...,b_{\sig(i)},\\
Q'_{j+k}(b_{\sig(i+1)},...b_{\sig(n)},m_{\tau(1)},...,m_{\tau(k)}),
m_{\tau(k+1)},...,m_{\tau(r)})\\[2mm]
+\sum_{k+l=r+1}\sum_{\sig\in\Sh(k,n)}\epsilon(\sig,b)
Q'_{n+r}(Q_k^B(b_{\sig(1)},...,b_{\sig(k)}),
b_{\sig(k+1)},...,b_{\sig(n)},m_1,...,m_r)
\end{spliteqn*}

But this term just equals
\begin{spliteqn*}
\sum_{k+l=n+r+1\atop k\geq n}
\sum_{\sig\in\Sh(k,n+r)\atop\sig(1)=1,...,\sig(n)=n}
\epsilon(\sig,u)
Q^M_l(Q'_k(b_1,...,b_n,u_{\sig(n+1)},...,u_{\sig(k)}),
u_{\sig(k+1)},...,u_{\sig(n+r)})\\
+\sum_{k+l=n+r+1}\sum_{\sig\in\Sh(k,n+r)\atop\sig(1)>n}
\epsilon(\sig,u)
Q'_l(Q^M_k(u_{\sig(1)},...,u_{\sig(k)}),b_1,...,b_n,
u_{\sig(k+n+1)},...,u_{\sig(n+r)})\\
+\sum_{k+l=n+r+1}\sum_{\sig\in\Sh(k,n+r)\atop\sig(k)\leq n}
\epsilon(\sig,u)
Q'_l(Q^B_k(b_1,...,b_k),b_{\sig(k+1)},...,b_{\sig(n)},m_1,...,m_r)\\[2mm]
+\sum_{k+l=n+r+1}
\sum_{\sig\in\Sh(k,n+r)\atop\sig(1)\leq n,\sig(k+1)\leq n}
\epsilon(\sig,u)
Q'_l(Q'_k(u_{\sig(1)},...,u_{\sig(k)}),u_{\sig(k+1)},...,u_{\sig(n+r)}),
\end{spliteqn*} which is
\begin{equation*}
(Q^M+Q^B+Q')^2(b_1,...,b_n,m_1,...,m_r).
\end{equation*}
Here, we have set $(u_1,...,u_{n+r}):=(b_1,...,b_n,m_1,...,m_r)$. 
This proves the second part of the theorem.\\

To prove its first part, we will
show that the map $F:=(f\times\id):B\times M\nach U\times M$
respects the DG structures, i.e. that
for $n\geq 0$ the following equality holds:

\begin{spliteqn}\label{star}
\tQ_1F_n+\frac{1}{2}\sum_{i+j=n}\tQ_2\circ F_i\ot F_j\circ\alpha_{i,n}
+\sum_{k=3}^n\sum_{I\in\NN^k\atop |I|=n}\frac{1}{I!k!}
\tQ_k\circ(F_{i_1}\ot...\ot F_{i_k})\circ\alpha_n=\\
\sum_{k+l=n+1}F_l\circ(\tQ'_k\ot 1\ot...\ot 1)\circ\alpha_{k,n}\nonumber
\end{spliteqn}

Remark that $F_n$ takes the following values on
products $b_1\odos b_r\odot m_1\odos m_{n-r}$, with $r<n$,
$b_i\in B$ and $m_j\in M$:
\begin{align*}
F_n(b_1,...,b_r,m_1,...,m_{n-r})=0 &\quad\quad \text{ for} & 0<r<n\\
F_n(b_1,...,b_n)=f_n(b_1,...,b_n) &&\\
F_n(m_1,...,m_n)=0 & \quad\quad\text{ for}& m>1\\
F_1(m_1)=m_1&&
\end{align*}

Applying the left hand-side of equation (\ref{star}) on 
$b_1\odos b_r\odot m_1\odos m_{n-r}$, we only get the term
$$Q_{1+n-r}(f_r(b_1,...,b_r),m_1,...,m_r).$$

Applying the right hand-side of equation (\ref{star}) on 
$b_1\odos b_r\odot m_1\odos m_{n-r}$, we only get the term
$$Q'_n(b_1,...,b_r,m_1,...,m_r).$$
By our construction, both terms coincide. 
\qed

At the end of Section \ref{hodgecomp} we will be able to
construct a semiuniversal deformation of an $L_\infty$-algebra
with split tangent complex (see Theorem~\ref{semiuni}).

\section{Trees}\label{trees}

Trees were used by Gugeheim/Stasheff \cite{GuSt}, Merkulov
\cite{MerkK}, Kontsevich/Soibelman \cite{KontS} and
others to construct infinity structures. 
We define binary trees (in a slightly different manner
as usual) and assign several invariants to them, which
are important to get good signs, later.

\subsection{Definitions}
\begin{defi}
A \textbf{tree} with $n$ leaves is a pair $\phi=(\phi,V)$
consisting of a set
$V=\{K_0\ddd K_{n-2} \}$ of \textbf{ramifications} such
that for each $i=0\ddd n-2$, we have:
\begin{enumerate}
\item
$\phi^{-1}(K_i)$ contains at most 2 elements.
\item
There is an $n\geq 0$ such that $\phi^n(K_i)=K_0$.
\end{enumerate}
$K_0$ is called \textbf{root} of $\phi$. 
\end{defi}

There is a tree with one leaf and no ramification,
which will always be denoted by $\tau$.

\begin{defi}
An \textbf{orientation} of a tree $(\phi,V)$ is a family
$\pi=(\pi_K)_{K\in V}$ of inclusions $\pi_K:\phi^{-1}(K)\nach
\{1,2\}$.
The triple $\phi=(\phi,V,\pi)$ is called an \textbf{oriented tree.}
\end{defi}

\begin{bem}
For each oriented tree $(\phi,V,\pi)$, there is a natural ordering
on the set $V$: For $K\in V\setminus K_0$, suppose that $\phi^m(K)=K_0$.
We set
$$v(K):=\frac{\pi_{\phi(K)}(K)}{3^m}+
\frac{\pi_{\phi^2(K)}(\phi(K))}{3^{m-1}}\kkk
\frac{\pi_{\phi^m(K)}(\phi^{m-1}(K))}{3}.$$
Set $v(K_0):=0$. Then $v:V\nach\RR$ is invective, hence it induces
an ordering on $V$.
\end{bem}

When we write down the value $v(K)$ of a ramification $K$ in its
3-ary decomposition, we just get an algorithm, how to get from the root $K_0$ 
to $K$. For example $0.1121$ means 
``go (in the driving direction) right-right-left-right''.
When $(\phi,V,\pi)$ is an oriented tree with $n$ leaves, we can extend
the map $\phi$ to a map 
$\tilde{\phi}:V\setminus{K_0}\cup\{1,...,n\}\nach V$, such
that
\begin{itemize}
\item
For $1\leq i< j\leq n$ we have $\tilde{\phi}(i)\leq\tilde{\phi}(j)$.
\item
For each $K\in V$, $\tilde{\phi}^{-1}(K)$ has exactly 2 elements.
\end{itemize}
 
The numbers $1,...,n$ stand for the leaves of $\phi$.
Furthermore, we can extend the map $v:K\nach (0,1)$ on
$\tilde{V}:=V\cup\{1,...,n\}$ in such a way that the 3-ary
decomposition of $v(i)$ describes the way from the root to the
i-th leaf of $\phi$, for $i=1,...,n$. Then we have
$v(i)<v(j)$ for $1\leq i<j\leq n$. In consequence, we have
an ordering on $\tilde{V}$.

\begin{defi}
Two trees $(\phi,V)$ and $(\phi',V')$ are called \textbf{equivalent}
if there is a bijection $f:V\nach V'$ of the ramification sets
such that $f\circ\phi=\phi'\circ f$.
Two oriented trees $(\phi,V,\pi)$ and $(\phi',V',\pi')$ are called
\textbf{oriented equivalent}
if there is a bijection $f:V\nach V'$ of the ramification sets
such that $f\circ\phi=\phi'\circ f$ and $\pi'\circ f=\pi$.  
\end{defi}
When we draw oriented trees, we shall put elements $K'$ of 
$\phi^{-1}(K)$ down left of $K$ if $\pi_K(K')=1$ and down right 
of $K$ if $\pi_K(K')=2$.

\begin{beisp}\label{extree}
The following trees with three leaves are equivalent
but not oriented equivalent: 

{\footnotesize
\setlength{\unitlength}{0.4cm}
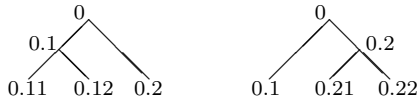
\begin{picture}(14,4)
\put(1,1){\line(1,1){2}}
\put(2,2){\line(1,-1){1}}
\put(3,3){\line(1,-1){2}}
\put(9,1){\line(1,1){2}}
\put(11,3){\line(1,-1){2}}
\put(11,1){\line(1,1){1}}
\put(0.3,0.5){$0.11$}
\put(2.5,0.5){$0.12$}
\put(1,2){$0.1$}
\put(2.5,3){$0$}
\put(8.4,0.5){$0.1$}
\put(4.5,0.5){$0.2$} 
\put(10.5,3){$0$}
\put(10.5,0.5){$0.21$}
\put(12.2,2){$0.2$}
\put(12.6,0.5){$0.22$}
\end{picture}
}

For each ramification and each leaf we have indicated its value.
\end{beisp}

Set $\Ot(n)$ to be the set of equivalence classes of oriented trees with
$n$ leaves.

\begin{beisp}
\begin{enumerate}
\item
$\Ot(2)$ contains just one element. It will always be
denoted by $\beta$.
\item
$\Ot(4)$ contains just the following elements: 

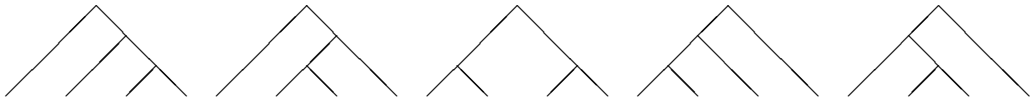
\setlength{\unitlength}{0.4cm}
\begin{picture}(39,5)
\put(1,1){\line(1,1){3}}
\put(3,1){\line(1,1){2}}
\put(5,1){\line(1,1){1}}
\put(4,4){\line(1,-1){3}}

\put(8,1){\line(1,1){3}}
\put(10,1){\line(1,1){2}}
\put(11,2){\line(1,-1){1}}
\put(11,4){\line(1,-1){3}}

\put(15,1){\line(1,1){3}}
\put(19,1){\line(1,1){1}}
\put(16,2){\line(1,-1){1}}
\put(18,4){\line(1,-1){3}}

\put(22,1){\line(1,1){3}}
\put(23,2){\line(1,-1){1}}
\put(24,3){\line(1,-1){2}}
\put(25,4){\line(1,-1){3}}

\put(29,1){\line(1,1){3}}
\put(31,1){\line(1,1){1}}
\put(31,3){\line(1,-1){2}}
\put(32,4){\line(1,-1){3}}
\end{picture}
\end{enumerate}
\end{beisp}

\begin{bem}
For a tree $(\phi,V)$ and $K\in V$, there is a tree $\phi|_K$
with root $K$ and ramifications
$\{ K'\in V:\;\phi^n(K')=K \text{ for an } n\geq 0\}$.
\end{bem}

We have to introduce several invariants:\\

For a tree $\phi$ with $n>1$ leaves and $1\leq i\leq n$, set $w_\phi(i)$
to be the difference of the number $s_\phi(i)$ of ramifications
of $\phi$ which are smaller than $i$ and $i-1$. ($i-1$ is the
number of leaves of $\phi$, smaller than $i$.) 
For $K\in V$ set $w_\phi(K):=w_{\phi-\phi|_K}(K)$, where on the right
hand-side $K$ is considered as leaf of $\phi-\phi|_K$.\\

\begin{bem}
For $K\in V$, the integer $w_\phi(K)$ is just the number of 1's
arising in the 3-ary decomposition of $v(K)$.
\end{bem}

Now, for each tree $\phi$ with at least 2 leaves,
set $e(\phi):=(-1)^{w_\phi(1)+...+w_\phi(n)}$.
Set $e(\tau):=1$

\begin{beisp}
\begin{enumerate}
\item
$e(\beta)=-1$
\item
For the first tree in Example~\ref{extree}, we have $e(\phi)=-1$;
For the second tree in Example~\ref{extree}, we have $e(\phi)=+1$;
\end{enumerate}
\end{beisp}

Now, let $L$ be a graded module, $\phi$ an oriented tree with $n$
leaves and $B=(b_K)_{K\in V}$ a family of bilinear maps $L\ot L\nach L$.
Recursively, we want to define a multilinear map
$$\phi(B):L^{\ot n}\nach L.$$

\begin{itemize}
\item
If $\phi$ has one leaf, i.e. $B$ is empty, we set $\phi(B):=\id$.
\item
If $\phi$ has only two leaves, i.e. $V=\{K_0\}$,
for a bilinear map $b_0: L\ot L\nach L$, we set
$\phi(b_0):=b_0$.
\item
If $\phi^{-1}(K_0)$ contains exactly one element, say $K_1$, and
$\pi_{K_0}(K_1)=1$, we set
$$\phi(B):=b_0\circ(\phi|_{K_1}((b_K)_{K\in V\setminus K_0})\ot 1).$$
\item
If $\phi^{-1}(K_0)$ contains exactly one element, say $K_1$, and
$\pi_{K_0}(K_1)=2$, we set
$$\phi(B):=b_0\circ(1\ot\phi|_{K_1}((b_K)_{K\in V\setminus K_0})).$$
\item
If $\phi^{-1}(K_0)=\{K_1,K_2\}$ with $\phi_{K_0}(K_1)=1$ and 
$\phi_{K_0}(K_2)=2$, we set
$$\phi(B):=b_0\circ
(\phi|_{K_1}((b_K)_{K\in V_1})\ot\phi|_{K_2}((b_K)_{K\in V_2})).$$
Here, $V_1$ denotes the ramification set of $\phi|_{K_1}$ and
$V_2$ the ramification set of $\phi|_{K_2}$.
\end{itemize}

\subsection{Operations on trees}

\paragraph{Addition}
Let $(\phi,V,\pi)$ and $(\phi',V',\pi')$ be oriented trees with 
disjoint ramification sets. Let $R$ be a point in neither one 
of them.
Set $V'':=V\cup V'\cup \{R\}$. We define a map 
$\psi:V''\setminus R\nach V''$
by $\psi|_{V\setminus K_0}:=\phi$, $\psi|_{V'\setminus K'_0}:=\phi'$
and $\psi(K_0):=\psi(K'_0):=R$.\\

There is a family $(\pi''_K)_{K\in W}$ of inclusions
$\pi''_K:\psi^{-1}(K)\nach\{0,1\}$ with
$\pi''_K=\pi_K$ for $K\in V$, $\pi''_K=\pi'_K$ for $K\in V'$ and
$\pi''_R(K_0)=0$ and $\pi''_R(K'_0)=1$.
Now, we set
$$(\phi,V,\pi)+(\phi',V',\pi'):=(\psi,V'',\pi'').$$
It is obvious, how to define the addition of non-oriented trees.
The addition of oriented trees is not commutative. The addition
of non-oriented trees is commutative.

\begin{beisp}
$\tau+\tau=\beta$. Furthermore, each tree can be reconstructed
by addition out of copies of $\tau$.
\end{beisp}

\paragraph{Subtraction}
Let $(\phi,V)$ be a tree with $n$ leaves and $K\in V$.
Let $l$ be the number of leaves of $\phi|_K$. Then the definition
of a tree $\phi-\phi|_K$ with $n-l+1$ leaves is quite obvious.

\paragraph{Composition}
Let $(\phi,V,\pi)$ be an oriented tree with $n$ leaves and
let\\
$(\psi^{(1)},V^{(1)},\pi^{(1)})\ddd (\psi^{(n)},V^{(n)},\pi^{(n)})$
be oriented trees. Let $W$ be the disjoint union of $V$ and all 
$V^{(i)}$.
For $K\in V$ set $n(K):=2-|\phi^{-1}(K)|$. (This is the number of leaves
belonging to $K$.) Let $K_1<\ldots<K_l$ all elements $K$ of $V$ with
$n(K)>0$. Then we define a map $\Phi:W\setminus K_0\nach W$ as follows:
For $K\in V\setminus K_0$ set $\Phi(K):=\phi(K)$.
For $K\in V^{(i)}\setminus K^{(i)}_0$ set $\Phi(K):=\psi^{(i)}(K)$.
And define the values of $\Phi$ on the $K^{(i)}$, setting
$$(\Phi(K^{(1)}_0)\ddd \Phi(K^{(n)}_0)):=
(\underbrace{K_1\ddd K_1}_{n_1\text{ times}},\ldots,
\underbrace{K_l\ddd K_l}_{n_l\text{ times}}).$$
Then, $(\Phi,W)$ is a tree with a canonical orientation $\pi'$, given as
follows:
For each $i$, $K\in V^{(i)}$ and $K'\in\Phi^{-1}(K)$, we set
$\pi'(K'):=\pi^{(i)}(K')$. For $K\in V$ and $K'\in\Phi^{-1}(K)\cap V$,
we set $\pi'(K'):=\pi(K')$. It remains to define $\pi'_{K_i}$ on elements
of $\Phi^{-1}(K_i)\setminus V$, for $i=1\ddd l$. So, if $n(K_i)$ equals
$2$, then $\Phi^{-1}(K_i)\setminus V$ has two elements, say
 $K^{(j)}_0$ and $K^{(k)}_0$ with $j<k$. Set $\Phi_{K_i}(K^{(j)}_0):=1$
and $\Phi_{K_i}(K^{(k)}_0):=2$.
If $n(K_i)$ equals
$1$, then $\Phi^{-1}(K_i)$ has one element in $V$, say $K$ and one 
element which is not in $V$, say $K'$.
Set $\Phi_{K_i}(K'):=1$ if $\phi_{K_i}(K)=2$ and
$\Phi_{K_i}(K'):=2$ if $\phi_{K_i}(K)=1$.\\

We will denote this decomposition by
$\Phi=\phi\circ(\psi^{(1)},...,\psi^{(n)})$.

\begin{bem}\label{exponent}
In this situation, suppose that there is a family $B=(b_K)_{K\in W}$ of
of bilinear maps $L\ot L\nach L$. Set $B^{(0)}:=(b_K)_{K\in V}$ and
$B^{(i)}:=(b_K)_{K\in V^{(i)}}$ for $i=1\ddd n$.
Then we have
$$\phi\circ(\psi^{(1)}\ddd\psi^{(n)})(B)=
(-1)^{\text{exponent} }
\phi(B^{(0)})\circ(\psi^{(1)}(B^{(1)})\ot\ldots \psi^{(n)}(B^{(i)})),$$
where the exponent is the sum
$(\sum_{K\in V^{(1)}}b_K)(\sum_{K\in V}^{V>1}b_K)+...+
\sum_{K\in V^{(n-1)}}b_K)(\sum_{K\in V}^{V>n-1}b_K)$.
We remind that $V>i$ means that the value $v(V)$ is greater
than the value $v(i)$ of the $i$-th leaf of $\phi$.
\end{bem}

\section{$L_\infty$-equivalence of $L$ and $H(L)$}\label{leq}

Let $L=(L,d,\lie)$ be a differential graded Lie-algebra, where the 
differential $d$ is of degree $+1$. Suppose that there is a splitting 
$\eta$, i.e. a map of degree $-1$ such that $d\eta d=d$. Furthermore, 
suppose that $\eta^2=0$ and $\eta d\eta=\eta$. When we use a Lie 
bracket on $\Hom(L,L)$, we mean the graded commutator. 

In this section, we want to construct an $L_\infty$-algebra structure
$\mu_\ast$ on $H:=H(L,d)$ with $\mu_1=0$,
such that $(L,d,[\cdot,\cdot])$ and $(H,\mu_\ast)$ are
$L_\infty$-equivalent.
The multilinear forms $\mu_n$ will be constructed using trees
as in the last section.
A similar construction for $A_\infty$-algebras 
can be found in \cite{GuSt} and \cite{MerkK}. 

We have to make some preparations. First of all, there is
the following easy but important remark:

\begin{bem}\label{sechszudrei}
Let $n\geq 3$ be a natural number. There is a 1:1-correspondence
between triples $(\Phi,K,\sigma)$, where $\Phi=(\Phi,V,\pi)$
is an oriented tree with $n$ leaves, K is a ramification in $V$,
$\sig$ a permutation in $\Sigma_n$ and 6-tuples
$(k,\phi,\psi,\rho,\gamma,\delta)$, where $k$ is a natural number
with $2\leq k\leq n-1$, $\phi$ is a tree in $\Ot(k)$,
$\psi$ is a tree in $\Ot(l)$ where $l:=n+1-k$, $\rho$ is a
shuffle in $\Sh(k,n)$ and $\gamma\in\Sigma_l$, $\delta\in\Sigma_k$
are permutations.\\

{\footnotesize
\setlength{\unitlength}{0.5cm}
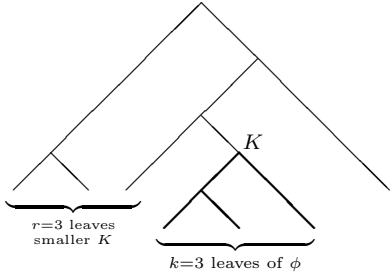
\begin{picture}(12,9)

\put(7.1,4.1){$K$}
\put(0.8,2.8){$\underbrace{\quad\quad\quad\quad\quad\quad}_{r=3
\text{ leaves} \atop \text{smaller }K}$} 
\put(4.8,1.8){$\underbrace{\quad\quad\quad\quad\quad\quad\quad}_{k=3
\text{ leaves of }\phi}$} 
\put(1,3){\line(1,1){5}}
\put(2,4){\line(1,-1){1}}
\put(4,3){\line(1,1){3.5}}
\put(6,8){\line(1,-1){5}}
\put(6,5){\line(1,-1){1}}
\thicklines
\put(5,2){\line(1,1){2}}
\put(6,3){\line(1,-1){1}}
\put(7,4){\line(1,-1){2}}
\put(13,5){\parbox{8cm}{\textbf{Example:}\\ 
The fine lines represent the tree $\psi$
and the fat lines the tree $\phi$.\\ 
In the sequel, the first $r$ leaves of $\Phi$ will be
associated to the indexes $\sig(1)=\rho(\gamma(1)+k-1),...
,\sig(r)=\rho(\gamma(r)+k-1)$, the following $k$ leaves
to the indexes $\sig(r+1)=\rho(\delta(1)),...,
\sig(r+k)=\rho(\delta(k))$ and the remaining leaves 
to the indexes $\sig(r+k+1)=\rho(\gamma(r+2)+k-1),...,
\sig(n)=\rho(\gamma(l)+k-1)$.}} 

\end{picture}}

Here-fore, to the triple $(\Phi,K,\sig)$, we associate the following
data:
Set $k$ to be the number of leaves of $\Phi|_K$, $\phi:=\Phi|_K$,
$\psi:=\Phi-\phi$. Let $r$ be the number of leaves $F$ of $\Phi$
with $F<K$. $\rho$ is chosen in such a way that
$\{\rho(1),...,\rho(k)\}=\{\sig(r+1),...,\sig(r+k)\}$.
$\delta$ is defined by $\delta(i):=\rho^{-1}(\sig(r+i))$ for
$i=1,...,k$ and $\gamma$ is defined in the following way:

$$\gamma(i):=\left\{ \begin{array}{r@{\quad\text{ for }\quad}l}
\rho^{-1}(\sig(i))-k+1 & i=1,...,r\\
1 & i=r+1\\
\rho^{-1}(\sig(i+k-1))-k+1 & i=r+2,...,l
\end{array}\right. $$ 

In the other way, to the 6-tuple $(k,\phi,\psi,\rho,\gamma,\delta)$,
we associate the following data:
Set $r:=\gamma^{-1}(1)-1$. Then $\Phi$ is the composition
$$\Phi=\psi\circ(\underbrace{\tau,...,\tau}_{r\text{ times}},\phi,
\tau,...,\tau),$$
where $\tau$ again stands for the tree with one leaf.
$K$ is the root of $\phi$, considered as ramification of
$\Phi$ and $\sig$ is given by

$$\sig(i):=\left\{ \begin{array}{r@{\quad\text{ for }\quad}l}
\rho(\gamma(i)+k-1) & i=1,...,r\\
\rho(\delta(i-r)) & i=r+1,...,r+k\\
\rho(\gamma(i-(k-1))+k-1) & i=r+k+1,...,n.
\end{array}\right. $$ 
\end{bem}

Now suppose that such corresponding tuples $(\Phi,\hat{K},\sig)$
and $(k,\phi,\psi,\rho,\gamma,\delta)$ are given.
Let $V'$ be the ramification set of $\psi$ and $V''$ the ramification
set of $\phi$. Then $V:=V'\cup V''$ is the ramification set of
$\Phi$. Again, set $r:=\gamma^{-1}(1)-1$. 
Remark that the ordering on $V$ depends on $\gamma$.
We define a permutation $\tilde{\gamma}\in\Sigma_{l-1}$ by
  
$$\tilde{\gamma}(i):=\left\{ \begin{array}{r@{\quad\text{ for }\quad}l}
\gamma(i)-1 & i=1,..., r\\
\gamma(i+1)-1 & i=r+1,...,l-1.
\end{array}\right. $$ 

\begin{lemma}\label{wiwo}
We keep all notations from above. Let
$B=(b_K)_{K\in V}$ be a family of homogeneous bilinear forms
$L\ot L\nach L$. Denote the subfamilies $(b_K)_{K\in V'}$ and
$(b_K)_{K\in V''}$ by $B'$ and $B''$. Set $W$ to be the set
of all ramifications $K\in V$ such that $K>\hat{K}$. Then we have
\begin{multline*}
\psi(B')\circ\gamma\circ(\phi(B'')\circ\delta\ot 1\ot...\ot 1)
\circ\rho\\[3mm]
=(-1)^{r+rk}\psi(B')\circ(\underbrace{1\ot...\ot 1}_{r\text{ times}}
\ot\phi(B'')\ot 1\ot...\ot 1)\circ\sig\\
=(-1)^{r+rk+\sum_{K\in W}b_k\cdot B''}\Phi(B)\circ\sig.
\end{multline*}
\end{lemma}
\bew
Let $a_1,...,a_n$
be homogeneous elements of $L$. 
Then we have
\begin{multline*}
(\psi(B')\circ\gamma\circ(\phi(B'')\circ\delta\ot 1\ot...\ot 1)
\circ\rho)(a_1\ot...\ot a_n)=\\[3mm]
=\chi(\rho,a_1,...,a_n)\chi(\delta,a_{\rho(1)},...,a_{\rho(k)})
\cdot(\psi(B')\circ\gamma)(
\underbrace{\phi(B'')(a_{\rho(\delta(1))},...,a_{\rho(\delta(k))})}_{=:u_1}
\ot\underbrace{a_{\rho(k+1)}}_{u_2}\ot...\ot
\underbrace{a_{\rho(n)}}_{u_l})\\
=\chi(\rho,a_1,...,a_n)\chi(\delta,a_{\rho(1)},...,a_{\rho(k)})
\chi(\gamma,u_1,...,u_l)\psi(B')(u_{\gamma(1)},...,u_{\gamma(l)}).
\end{multline*}
Using the following three formulas

\begin{align*}
\chi(\sig,a_1,...,a_n)&=(-1)^{kr+(a_{\sig(1)}+...+a_{\sig(r)})
(a_{\rho(1)}+...+a_{\rho(k)})}\chi(\rho,a_1,...,a_n)\cdot\\
&\quad\hspace{3cm}\chi(\delta,a_{\rho(1)},...,a_{\rho(k)})
\chi(\tilde{\gamma},a_{\rho(k+1)},...,a_{\rho(n)}),\\[2mm]
u_{\gamma(1)}\ot...\ot u_{\gamma(l)}&=
(-1)^{B''(a_{\sig(1)}+...+a_{\sig(r)})}
(1\ot...\ot 1\ot\phi(B'')\ot 1...\ot 1)
(a_{\sig(1)}\ot...\ot a_{\sig(n)}),\\[2mm]
\chi(\gamma,u_1,...,u_l)&=
(-1)^{r+u_1(u_{\gamma(1)}+...+ u_{\gamma(r)})}
\chi(\tilde{\gamma},u_2,...,u_l),
\end{align*}
this expression is just
$$(-1)^{kr+r}\chi(\sig,a_1,...,a_n)\psi(B')
((1\ot...\ot 1\ot\phi(B'')\ot 1\ot...\ot 1)
(a_{\sig(1)}\ot...\ot a_{\sig(n)})).$$
The second equality of this Lemma is just a special case of
Remark~\ref{exponent}.
\qed

We turn to the construction of an $L_\infty$-structure
on $H(L)$.

\begin{bem}
$[d,\eta]=d\eta+\eta d$ is a projection, i.e. $[d,\eta]^2=[d,\eta]$. 
And $H:=\Kern[d,\eta]$ is, as module, isomorphic to $H(L)$. The bracket 
on $L$ induces a Lie-bracket on $H(L)$ and the induced bracket on $H$ 
(via the isomorphism $H\nach H(L)$) is just given by 
$(1-d\eta)\lie=(1-[d,\eta])\lie$.
\end{bem}

For simplicity, we set $g:=\eta\lie$.

\begin{satz}~\label{trick}
The following graded anti-symmetric maps $\mu_n:H^{\ot n}\nach H$ 
of degree $2-n$ define
the structure of an $L_\infty$-algebra on $H$:
\begin{align*}
\mu_1:&=0\\
\mu_2:&=(1-d\eta)\lie\\
&\vdots\\
\mu_n:&=(\frac{-1}{2})^{n-1}\sum_{\phi\in\Ot_n}e(\phi)\phi((1-[d,\eta])\lie,g
\ddd g)\circ\alpha_n
\end{align*} 
\end{satz}
\bew
We must show that
\begin{equation}
\sum_{k+l=n+1}(-1)^{k(l-1)}\mu_l\circ(\mu_k\ot1\ot...\ot 1)
\circ\alpha_{k,n}=0.
\end{equation}
Up to the factor $(-1)^{n-1}$ this sum has the form
\begin{multline}\label{nonc}
\sum_k
\sum_{\phi,\psi}\sum_{\rho,\gamma,\delta}(-1)^{k(l-1)}e(\phi)e(\psi)
\psi((1-[d,\eta])\lie,g,...,g)\circ\gamma\circ\\
\circ
(\phi((1-[d,\eta])\lie,g,...,g)\circ\delta\ot1\ot...\ot1)\circ\rho,
\end{multline}
where $k$ ranges from from $2$ to $n-1$, $l=n+1-k$,
$\phi$ and $\psi$ vary in $\Ot(k)$ and $\Ot(l)$, $\rho$ in
$\Sh(k,n)$, $\gamma$ and $\delta$ in $\Sigma_l$ and $\Sigma_k$.
For corresponding tuples
$(k,\phi,\psi,\rho,\gamma,\delta)$ and $(\Phi,\hat{K},\sig)$
as in Remark~\ref{sechszudrei}, we denote
as usual $r:=\gamma^{-1}(1)-1$ and by $t$ the number of ramifications
of $\psi$, greater than $r+1$.
Using
\begin{eqnarray*}
e(\Phi)&=&(-1)^{w_{\Phi}(\hat{K})(k-1)}e(\phi)e(\psi)\\
w_{\Phi}(\hat{K})&=&l -1-r-t
\end{eqnarray*} 
and Lemma~\ref{wiwo}, the expression~\ref{nonc} can be expressed as
$$\sum_{\Phi\in\Ot(n)}\sum_{\hat{K}\in V\setminus K_0}e(\Phi)
(-1)^{r+w_\Phi(\hat{K})}\Phi(B)\circ\alpha_n,$$
where $B=(B_K)_{K\in V}$ is the family with 
$b_{K_0}=b_{\hat{K}}=(1-[d,\eta])\lie$ and $b_k=\eta\lie$
for $K\neq K_0,\hat{K}$.
To show that the last term is zero, it is enough to show
the following two conditions:
\begin{itemize}
\item
$\sum_{\Phi\in\Ot(n)}\sum_{K\in V\setminus K_0}
\sum_{\sig\in\Sigma_n} (-1)^{r+w_\Phi(K)}e(\Phi)
\Phi((1-[d,\eta])\lie,g,...,g,\underbrace{\lie}_{\text{pos. }K},
g,...,g)\circ\sigma=0.$
\item
For each tree $\Phi$, we have
\begin{equation}\label{nomal}
\sum_{K\in V\setminus K_0}
(-1)^{r+w_\Phi(K)} e(\Phi)
\Phi((1-[d,\eta])\lie,g,...,g,\underbrace{[d,\eta]\lie}_{\text{position }K},
g,...,g)\circ\sigma=0.
\end{equation}
\end{itemize}
The first condition follows by the Jacobi-identity and an easy
combinatorial argument. In equation (\ref{nomal}) the term
annihilate each other since the differential d trickles down
the branches of $\Phi$:\\ 

\textbf{Initiation of the trickling:}
Suppose that $\Phi^{-1}(K_0)$ contains an element
$K'$ with $\pi_{K_0}(K')=1$. Then we have the following picture:
 
\setlength{\unitlength}{0.5cm}
{\footnotesize
\begin{picture}(23,5)
\put(0,2){$0$}
\put(1,2){$=$}
\put(3,2){\line(1,1){1}}
\put(4,3){\line(1,-1){1}}
\put(3,3.2){$(1-[d,\eta])d\lie$}
\put(2,1.2){$\eta\lie$}
\put(5,1.2){$\eta\lie$}

\put(7,2){$=$}

\put(10,2){\line(1,1){1}}
\put(11,3){\line(1,-1){1}}
\put(10,3.2){$(1-[d,\eta])\lie$}
\put(9,1.2){$d\eta\lie$}
\put(12,1.2){$\eta\lie$}

\put(14,2){$+(-1)^{\sharp\text{ramific. of }\Phi|_{K'}}$}
 
\put(20,2){\line(1,1){1}}
\put(21,3){\line(1,-1){1}}
\put(20,3.2){$(1-[d,\eta])\lie$}
\put(19,1.2){$\eta\lie$}
\put(22,1.2){$d\eta\lie$}
\end{picture}
}

Here, we only have drawn the top of the tree $\Phi$
for the case where $\Phi^{-1}(K_0)$ consists of two
elements $K',K''$ and the corresponding bilinear
forms. It is quite evident how this goes  when
$\Phi^{-1}(K_0)$ has only one element, since
$d|_H=0$.\\

\textbf{Going-on of the trickling} at a ramification
$K\in V$: 

We illustrate the case, where $\Phi^{-1}(K)$ has two elements $K',K''$
with $\pi_K(K')=1$.

{\footnotesize
\begin{picture}(23,5)

\put(2,2){\line(1,1){1}}
\put(3,3){\line(1,-1){1}}
\put(2,3.2){$\eta d\lie$}
\put(1,1.2){$\eta\lie$}
\put(4,1.2){$\eta\lie$}

\put(6,2){$-$}

\put(9,2){\line(1,1){1}}
\put(10,3){\line(1,-1){1}}
\put(9.5,3.2){$\eta\lie$}
\put(8,1.2){$d\eta\lie$}
\put(11,1.2){$\eta\lie$}

\put(12,2){$-(-1)^{\sharp\text{ramific. of }\Phi|_{K'}}$}
 
\put(18,2){\line(1,1){1}}
\put(19,3){\line(1,-1){1}}
\put(18.5,3.2){$\eta\lie$}
\put(17,1.2){$\eta\lie$}
\put(20,1.2){$d\eta\lie$}

\put(21,2){$=\quad 0$}
\end{picture}
}

Iterating the trickling down to the leaves and
using $d|_H=0$, we see that all terms in the sum are annihilated.
\qed

\begin{satz}
The following anti-symmetric maps $f_n:H^{\ot n}\nach L$ of degree
$n-1$ define an $L_\infty$-equivalence $H\nach L$ (i.e. an 
$L_\infty$-quasi-isomorphism).
\begin{align*}
f_1:&=\text{inclusion}\\
f_2:&=-g\\
&\vdots\\
f_n:&=-(\frac{-1}{2})^{n-1}\sum_{\phi\in\Ot(n)}e(\phi)
\phi(g\ddd g)\circ\alpha_n.
\end{align*}
\end{satz}

\bew
For $n\geq 0$, we have to prove the equation
\begin{equation*}
df_n-\sum_{i+j=n}\frac{(-1)^i}{2}[f_i,f_j]\alpha_{i,n}=
\sum_{k+l=n+1}(-1)^{k(l-1)}f_l\circ(\mu_k\ot 1\ot...\ot 1)
\circ\alpha_{k,n}
\end{equation*} 
For $l=1$, the right hand-side is just $\mu_n$. Since
$$df_n=(\frac{-1}{2})^{n-1}\sum_{\phi\in\Ot(n)}
e(\phi)\phi(-d\eta\lie,g\ddd g)\circ\alpha_n,$$
it is sufficient to show the following three identities:

\begin{equation}\label{eins}
-\sum_{i+j=n}\frac{(-1)^{i}}{2}[f_i,f_j]\alpha_{i,n}=
(\frac{-1}{2})^{n-1}\sum_{\phi_\in\Ot(n)}e(\phi)
\phi(\lie,g\ddd g)\circ\alpha_n
\end{equation}
\begin{equation}\label{zwei}
f_l\circ(\phi(\lie,g\ddd g)\circ\alpha_k\ot 1\ot...\ot 1)
\circ\alpha_{k,n}=0 \text{ for }  l>1,k+l=n+1.
\end{equation}
\begin{multline}\label{drei}
(\frac{-1}{2})^{n-1}\sum_{\phi_\in\Ot(n)}e(\phi)
\phi(\eta d\lie,g\ddd g)\circ\alpha_n=
-\sum_{k+l=n+1 \atop l\geq 2}
(-1)^{k(l-1)}\sum_{\phi\in\Ot(k)}(\frac{-1}{2})^{k-1}e(\phi)\\
f_l\circ(\phi(d\eta+\eta d)\lie,g\ddd g)\circ\alpha_k\ot1\ot...\ot 1)
\circ\alpha_{k,n}.
\end{multline}

\textbf{Proof of equation (\ref{drei})}
The right hand-side of equation (\ref{drei}) is
$$(\frac{-1}{2})^{n-1}
\sum_{k+l=n+1}^{k,l\geq 2}\sum_{\phi,\psi}\sum_{\gamma,\delta,\rho}
(-1)^{k(l-1)}e(\phi)e(\psi)\psi(g,...,g)\circ\gamma\circ
(\phi([d,\eta]\lie,g,...,g)\circ\delta\ot 1\ot...\ot 1)\circ\rho.$$
As in the proof of Theorem~\ref{trick},
this expression takes the form
$$(\frac{-1}{2})^{n-1}
\sum_{\Phi\in\Ot(n)}\sum_{\hat{K}\in V}\sum_{\sig\in\Sigma_n}
(-1)^{r+w_\Phi(\hat{K})}e(\Phi)\Phi(B)\sig,$$
where $B=(b_K)_{K\in V}$ is
the family with $b_{\hat{K}}=[d,\eta]\lie$ and
$b_K=\eta\lie$ for $K\neq\hat{K}$.\\

Hence to show equation (\ref{drei}), it is enough to show that 
for each tree $\phi$, we have
$$\Phi(\eta d\lie,g,...,g)=\sum_{K\in V\setminus K_0}
(-1)^{r+w_\Phi(K)}\Phi(B).$$
This is true by the same trickling argument as in 
Theorem~\ref{trick}.\\

\textbf{Proof of equation (\ref{zwei}):}
This is again the Jacobi-identity and some combinatorics.\\

\textbf{Proof of equation (\ref{eins}):}  
\begin{multline*}
\sum_{i+j=n}\frac{(-1)^i}{2}\lie\circ(f_i\ot f_j)\circ\alpha_{i,n}=\\
=\sum_{i+j=n}(-\frac{1}{2})^{i-1+j-1}
\sum_{\phi\in\Ot(i),\psi\in\Ot(j)}\frac{(-1)^i}{2}
e(\phi)e(\psi)(\phi+\psi)(\lie,g,...,g)\circ\alpha_n\\
=-(-\frac{1}{2})^{n-1}\sum_{\Phi\in\Ot(n)}
e(\Phi)\Phi(\lie,g,...,g)\circ\alpha_n.
\end{multline*}
\qed

\section{Decomposition Theorem for Differential 
Graded Lie Algebras}\label{hodgecomp}

In this section, we want to realize two things:
(a) the construction of an inverse map of the quasi-isomorphism
$f:(H,\mu_\ast)\nach (L,d,\lie)$, constructed in Section~\ref{leq};
(b) the construction of a semi-universal deformation $(V,Q^V,Q')$
for a given formal DG manifold.
As consequence of (a), we get the following decomposition theorem
$(L,d,\lie)\isom(H,\mu_\ast)\oplus (F,d,0)$ for DGLs, where $F$
is the complement of $H$ in $L$ and the sum is taken in the 
category of $L_\infty$-algebras. The existence of such a decomposition,
was already stated by Kontsevich (see \cite{Kont}) and an $A_\infty$-analogue
was proved by Kadeishvili (see \cite{Kade}). In fact, each $L_\infty$-
(resp. $A_\infty$-algebra) over a field
is isomorphic to the direct sum of a 
minimal and a linear contractible one. 
Our Proposition~\ref{quill} is analogue to  
the corresponding statement for $A_\infty$-algebras, which
was proved by Lefevre (see \cite{Lefe}). The proof here is almost
a transcription of Lefevre's proof.

\subsection{Obstructions}\label{obstruct}
Consider the formal DG manifolds  $(W,Q)$ and $(W',Q')$. For
any $n\geq 0$, there is a differential $\delta$ of degree $+1$
on the graded module 
$\Hom(W^{\ot n},W')$, given
by $\delta(g)=Q_1'\circ g-(-1)^gg\circ Q_1^{(n)}$.
Now, let $f:W\nach W'$ be a morphism of formal supermanifolds.
Set
$$r(f_1,...,f_{n-1}):=\sum_{k+l=n+1\atop k\geq 2}f_l\circ Q_k^{(n)}
-\sum_{k=2}^n\sum_{i_1+...+i_k=n}Q'_k\circ f_I.$$

Recall that $f$ is an $L_\infty$-homomorphism iff, for each $n\geq 1$,
we have $\delta(f_n)=r(f_1,...,f_n)$. If this condition is
satisfied only for $n\leq m$, we will call $f$ (or the family
$(f_1,...,f_m)$) an \textbf{$L_m$-homomorphism}.

\begin{lemma}\label{obst}
Suppose that $f$ is an $L_{n-1}$-homomorphism. Then
$\delta(r(f_1,...,f_{n-1}))=0$.
\end{lemma}

The proof is done in the appendix.

\begin{bem}\label{obem}
Let $e:W\nach W'$ and $f:V'\nach V$ be strict $L_\infty$-morphisms 
and let $g:V\nach W$ be any $L_\infty$-morphism.
Then
\begin{enumerate}
\item
$r((gf)_1,...,(gf)_{n-1})=r(g_1,...,g_{n-1})\circ f_1^{\ot n}.$
\item
$r((eg)_1,...,(eg)_{n-1})=e_1\circ r(g_1,...,g_{n-1})$.
\end{enumerate}
\end{bem}

\subsection{Constructions}

\begin{prop}\label{isom}
Let $f:M\nach M'$ be a morphism of formal supermanifolds.
Suppose, there is a module homomorphism $g':M'\nach M$ such that
$g'\circ f_1=\id_M$. Then, there is a morphism $g:M'\nach M$
of formal supermanifolds such that $g_1=g'$ and $gf=\id_M$.
If $f_1$ is an isomorphism with inverse $g$ and if
$f$ is $Q$-equivariant and if
$g'$ respects $Q_1^{M'}$ and $Q_1^M$, then $g$ can be chosen 
$Q$-equivariant, as well.
\end{prop}
\bew
One can check directly that the 
sequence of maps defined by
$$g_n:=-\sum_{k=2}^n\sum_{I\in\NN^k\atop |I|=n}g_1\circ f_k\circ (g_{I}),$$
for $n\geq 2$, define a morphism of formal supermanifolds
with the desired property.
\qed

\begin{lemma}\label{strick}
Let $f:V\nach W$ be a morphism of formal DG manifolds.
\begin{enumerate}
\item
If $f_1$ is split injective, then there is a formal DG manifold
$W'$ and an $L_\infty$-isomorphism $\kappa:W\nach W'$ such that
$\kappa\circ f$ is strict. 
\item
If $f_1$ is split surjective, then there is a formal DG 
manifold
$V'$ and an $L_\infty$-isomorphism $\kappa:V'\nach V$ such that
$f\circ\kappa$ is strict. 
\end{enumerate}
\end{lemma}
\bew
(i) As module, set $W':=W$.
We have to construct an isomorphism $\kappa:S(W)\nach S(W')$ of
graded coalgebras and then, we can define the DG structure on
$W'$ via 
$Q^{W'}:=\kappa\circ Q^W\circ\kappa^{-1}$. 
Set $\kappa_1:=\id$. Inductively, we define
maps $\kappa_n:W^{\odot n}\nach W'$ such that for $2\leq m\leq n$, we
have 
$$(\kappa\circ f)_m=\sum_{k=1}^m\sum_{I\in \NN^k\atop |I|=n}
\kappa_k\circ f_I=0.$$
Let $g:W\nach V$ be a module homomorphism with
$g\circ f_1=\id_V$.
When $\kappa_1,...,\kappa_n$ is already constructed,
set
$$\kappa_{n+1}:=-\sum_{k=1}^n\sum_{I\in \NN^k\atop |I|=n+1}
\kappa_k\circ f_I\circ g^{\odot n+1}.$$
Obviously, $(\kappa\circ f)_m=0$ for $2\leq m\leq n+1$.
(ii) goes in a similar way.
\qed

For our situation, we have the following, more explicit
statement:

\begin{lemma}
Let $f:H\nach L$ be the $L_\infty$-quasi-isomorphism, constructed in
Section~\ref{leq}. 
Consider the applications $\kappa_n:L^{\ot n}\nach L$,
defined by $\kappa_1:=\id$ and $\kappa_n:=-f_n\circ\pr_H^{\ot n}$,
for $n\geq 2$.
Then, $\kappa$ is an $L_\infty$-morphism and
$\kappa\circ f$ is strict. Furthermore,
$(\kappa^{-1})_1=\id$ and $(\kappa^{-1})_n=f_n\circ\pr_H^{\ot n}$
for $n\geq 2$. 
\end{lemma}
\bew
For $n\geq 2$, we have
\begin{equation*}
(\kappa\circ f)_n=\sum_{k=1}^n\sum_{I\in\NN^K\atop|I|=n}\kappa_k\circ f_I
=f_n-f_n=0.
\end{equation*}
The second statement is as easy to prove.
\qed

The following important proposition says that the quadruple
(category of $L_\infty$-algebras; class of $L_\infty$-quasi-isomorphisms;
class of those $L_\infty$-morphisms $f$ such that $f_1$ is split
injective; class of those $L_\infty$-morphism $f$ such that $f_1$ is split
surjective) satisfies one of Quillen's axioms for model categories. 

\begin{prop}\label{quill}
Let 
$$\xymatrix{
A\ar[r]^c\ar[d]^f & C\ar[d]^e\\
B\ar[r]^d & D}$$
a commutative diagram of $L_\infty$-algebras.
Suppose that $f$ is split injective and that $e$ is split
surjective and that or $f$ or $e$ is an $L_\infty$-quasi-isomorphism.
Then, there is an $L_\infty$-morphism $g:B\nach C$, such that
the complete diagram
$$\xymatrix{
A\ar[r]\ar[d] & C\ar[d]\\
B\ar[r]\ar[ur]^g & D}$$
commutes.
\end{prop}
\bew 
By Lemma~\ref{strick}, 
we may suppose that $e$ and $f$ are strict.
Inductively, we will construct morphisms $g_n:B^{\odot n}\nach C$,
such that
\begin{enumerate}
\item
$\delta(g_m)+r(g_1,...,g_{m-1})=0$,
\item
$g_{m}\circ f_1^{\ot m}=c_m$,
\item
$e_1\circ g_m=d_m$,
\end{enumerate}
for each $m\leq n$.
Choose maps $u:(D,Q_1^D)\nach(C,Q_1^C)$ and
$v:(B,Q_1^B)\nach(A,Q_1^A)$ of DG-modules, such that
$v\circ f_1=\id_A$ and $e_1\circ u=\id_D$.
A candidate for $g_1$ can easily be found.
Suppose that $g_1,..,g_{n-1}$ are already constructed.
Then
$$\beta:=c_nv^{\odot n}+ud_n-ue_1c_nv^{\odot n}$$
satisfies conditions (i) and (ii).
By Lemma~\ref{obem}, we get
\begin{align*}
(\delta(\beta)+r(g_1,...,g_{n-1}))\circ f_1^{\odot n}=\\
\delta(\beta\circ f_1^{\odot n})+r((gf)_1,...,(gf)_{n-1})=\\
\delta(c_n)+r(c_1,...,c_{n-1})=0.
\end{align*}
On the other side, again by Lemma~\ref{obem}, we have
\begin{align*}
e_1\circ(\delta(\beta)+r(g_1,...,g_{n-1}))=\\
\delta(e_1\beta)+r((eg)_1,...,(eg)_{n-1})=\\
\delta(d_n)+r(d_1,...,d_{n-1})=0.
\end{align*}
Hence, $\delta(\beta)+r(g_1,...,g_{n-1})$ has
a factorization
$$\xymatrix{
B^{\ot n}\ar[r]^p &\Kokern(f_1^{\odot n})\ar[r]^ q &
\Kern(e_1)\ar[r]^i & C,
}$$
where $i$ is the natural inclusion and $p$ the natural epimorphism. 
By Lemma~\ref{obst}, $\delta(\beta)+r(g_1,...,g_{n-1})$
is a cycle, so $\delta(q)=0$, i.e. $q$
is a map of complexes.
Now, or $\Kokern(f_1^{\odot n})$ or $\Kern(e_1)$ is contractible.
Hence $q=\delta(h)$, for a morphism
$h:\Kokern(f_1^{\odot n})\nach\Kern(e_1)$ of graded modules.
Then $g_n:=\beta-i\circ h\circ p$ satisfies the conditions (i)-(iii).
\qed

\begin{kor}
There is a map $g:(L,d,\lie)\nach(H,\mu_\ast)$ of $L_\infty$-algebras
such that $g\circ f=\id_H$. 
\end{kor}

\begin{kor}\label{tride}
Let $M$ be an $L_\infty$-manifold and $(B,Q^B,Q)$ a deformation of
$M$, such that $(B,Q^B_1)$ is contractible and $Q_1=0$.
Then $(B,Q^B,Q)$ is a trivial deformation.
\end{kor}
\bew
There is a 
commutative diagram
$$\xymatrix{
M\ar[r]\ar[d] & (B\times M,Q^B+Q^M+Q)\ar[d]\\
(B\times M,Q^B+Q^M)\ar[r] & B
}$$
where the vertical left arrow induces an injective
quasi-isomorphism of DG-modules and the
vertical right arrow induces a surjective map of DG modules.
By Proposition~\ref{quill}, there is a map
$q:(B\times M,Q^B+Q^M)\nach (B\times M,Q^B+Q^M+Q)$ with $q_1=\id$,
completing the diagram commutatively.
In particular, $q$ establishes an isomorphism of the given deformation
and of the trivial deformation of $M$ with base $B$.
\qed

\begin{prop}
There exists 
\begin{enumerate}
\item
a homomorphism $\iota:(F,d,0)\nach (L,d,\lie)$
of $L_\infty$-algebras such that $\iota_1$ is the natural inclusion;
\item
a homomorphism $p:(L,d,\lie)\nach(F,d,0)$ of $L_\infty$-algebras
such that $p\circ\iota=\id_F$.
\end{enumerate}
\end{prop}
\bew
(a) Suppose, that there are already homomorphisms
$\iota_m:F^{\ot m}\nach L$, for $m\leq n-1$, which form
an $L_{n-1}$-homomorphism.
We have to find an $\iota_n$
such that $\delta(\iota_n)=r(\iota_1,...,\iota_{n-1})$.
Since $(F,d)$ is contractible, $\Hom(F^{\ot n},L)$ is
acyclic, so the existence of $\iota_n$ follows by
Lemma~\ref{obst}.\\
(b) By Lemma~\ref{strick}, we can assume that 
$\iota$ is strict. 
Set $p_1:=\pr_F=[d,\eta]$.
Now assume that $p_1,...,p_{n-1}$ are already constructed such
that they define an $L_{n-1}$-homomorphism $p':L\nach F$ such that
$(p'\circ\iota)_m=0$, for $m\leq {n-1}$. We have to find
$p_n:L^{\ot n}\nach F$ such that
\begin{align*}
\delta(p_n)+r(p_1,...,p_{n-1})&=0,\\
p_n\circ\iota^{\ot n}&=0.
\end{align*} 
We may chose $p_n:=\eta r(p_1,...,p_{n-1})$.
Then, since $r(p_1,...,p_{n-1})\in\Kern\delta$, we have
$\delta(p_n)=[\eta,d]\circ r(p_1,...,p_{n-1})=r$, and again
by Remark~\ref{obem}, we get
$$p_n\circ\iota^{\ot n}=\eta\circ r(p_1,...,p_{n-1})\circ\iota^{\ot n}=0.$$
So inductively, the map $p$ can be constructed.
\qed

As consequence, we get the expected decomposition theorem
for differential graded Lie algebras admitting a splitting: 

\begin{satz}\label{decomposition}
We have an isomorphism of $L_\infty$-algebras
$$f\times\iota:H\times F\nach L.$$
\end{satz}

\begin{kor}
If $(L,\mu_\ast)$ and $(L',\mu'_\ast)$ are $L_\infty$-algebras
such that $(L;\mu_1)$  and  $(L';\mu'_1)$ are split, then,
for each $L_\infty$-quasi-isomorphism $f:L\nach L'$, there
exists an $L_\infty$-morphism $g:L'\nach L$, such that
$f_1$ and $g_1$ are inverse maps on the homology. In particular, if
$k$ is a field, then $L_\infty$-quasi-isomorphism is an equivalence
relation.
\end{kor}

\subsection{A Semiuniversal Deformation}

\begin{bem}\label{popo}
\begin{enumerate}
\item
Let $M=(M,Q^M)$ be a formal DG manifold and $N$ a $Q^M$-closed
submodule of $M$, i.e. for all $j\geq 1$ and
$n_1,...,n_j\in N$, we have $Q^M_j(n_1,...,n_j)\in N$.
Then, $(N,Q^M|_N)$ is a formal DG manifold and the inclusion
$N\nach M$ is a morphism in $\dgmanf$.
\item
Let $(B,Q^B,Q)$ be a deformation of $(M,Q^M)$. Suppose that
$(B,Q^B)$ is a direct sum of formal DG manifolds
$(B',Q^{B'})$ and $(B'',Q^{B''})$. Then, $(B'',Q^{B''},Q|_{B''\times M})$
is also a deformation of $M$ and the canonical map
\begin{equation}\label{qqmb}
(B''\times M,Q^{B''}+Q^M+Q|_{B''\times M})\nach 
(B\times M,Q^B+Q^M+Q)
\end{equation}
defines a morphism of deformations.
\item
If in the situation of (ii), $(B',Q^{B'}_1)$ is contractible, then
the map (\ref{qqmb}) is an equivalence of deformations.
\end{enumerate}
\end{bem}
\bew
The statements (i) and (ii) are easy to see.
To show (iii), we apply Proposition~\ref{quill} to
the commutative diagram
$$\xymatrix{
(B''\times M,Q^{B''}+Q^M+Q|_{B''\times M})\ar[d]\ar[r] &
(B\times M,Q^B+Q)\ar[d]\\
(B'\times(B''\times M),Q^{B'}+Q^{B''}+Q^M+Q|_{B''\times M})\ar[r] & B''
}$$
of formal DG manifolds.
We get an isomorphism
$$(B,Q^B,Q|_{B''\times M})\nach (B,Q^B,Q)$$
of deformations with base $B$.
Obviously, the left one is equivalent to $(B'',Q^{B''},Q|_{B''\times M})$.
\qed

For the rest of this subsection, we work in the setting of
Section~\ref{defolinf}.
Thus $L$ is the DGL \\
$\Coder(S(M),S(M))$ for some formal DG manifold
$(M,Q^M)$. Again, we must assume that the complex
$(L,d)$ has a splitting $\eta$.
Equip the homology $H$ of $(L,d)$ with the $L_\infty$-structure
$\mu_\ast$, constructed in Section~\ref{leq}.
Set $U:=L[1]$, $V:=H[1]$ and
denote the morphism $V\nach U$, induced by
the quasi-isomorphism $H\nach L$, constructed in Section~\ref{leq}
again by $f$.\\

Let $(U,Q^U,Q)$ be the universal deformation of $M$
(see Section~\ref{defolinf}).
Again, set $\tilde{Q}:=Q^U+Q^M+Q$.
By base change $f:V\nach U$, we get a deformation 
$(V,Q^V,Q')$ of $(M,Q^M)$.
Explicitly, on products
$v_1\odos v_r\odot m_1\odos m_s$ with $r,s\geq 1$ and $n=r+s$, the
perturbation $Q'$ is given by
$$Q'_n(v_1,...,v_r,m_1,...,m_s)=(\pno f_r(v_1,...,v_r))_s(m_1,...,m_s).$$
Set $\tilde{Q}':=Q^V+Q^M+Q'$.

\begin{satz}\label {semiuni}
The deformation $(V,Q^V,Q')$ is semi-universal.
\end{satz}
\bew
Since $(H,\mu_\ast)$ is minimal and $(U,Q^U,Q)$ is universal,
we only have to show, that there exists a morphism
of deformations from $(U,Q^U,Q)$ to $(V,Q^V,Q')$.
This is a consequence of Theorem~\ref{decomposition}
and Remark~\ref{popo}.
\qed

\appendix

\section{Appendix}
\subsection{Some Calculations}

\textbf{Proof of Proposition~\ref{Fcoal}:} 
Induction on n\\
The case $n=1$ follows by the commutativity of diagram
\ref{qudiag}
and $\Kern(\Delta^+)=W'$. Now suppose that the formula is proved for
all $m\leq n-1$. Then we have
\begin{multline*}
( F\ot F\circ\Delta^+)(w_1\ddd w_n)=\\[3mm]
\sum_{j=1}^{n-1}\sum_{\tau\in\Sh(j,n)}\epsilon(\tau,w_1\ddd w_n)
 \hat{F}_j(w_{\tau(1)}\ddd w_{\tau(j)})\ot \hat{F}_{n-j}
(w_{\tau(j+1)}\ddd w_{\tau(n)})=\\
\sum_{j=1}^{n-1}\sum_{k,k'}\sum_{I,I'}\frac{1}{k!k'!I!I'!}
[(F_{i_1}\odos F_{i_k})\circ\alpha_j\ot(F_{i'_1}\odos F_{i'_{k'}})
\circ\alpha_{n-j}]
\circ\alpha_{j,n}(w_1\odos w_n),
\end{multline*}
where $k$ ranges over $1,...,j$; $k'$ over $1,...,n-j$, $I$ takes all
values in $\NN^{k}$, such that $|I|=j$ and
$I'$ takes all values in $\NN^{k'}$, such that $|I'|=n-j$.
The last expression equals
\begin{equation}
\sum_{k=2}^n\sum_{l=1}^{k-1}\sum_{I\in\NN^k}^{|I|=n}
\frac{1}{I!l!(k-l)!}[(F_{i_1}\odos F_{i_l})\ot
(F_{i_{l+1}}\odos F_{i_k})]\circ\alpha_n(w_1,...,w_n)
\end{equation}

On the other side, we have
\begin{multline*}
\Delta^+(\sum_{k=2}^n\sum_{I\in\NN^k\atop |I|=n}\sum_{\sig\in\Sigma_n}
\frac{1}{I!k!}\epsilon(\sig)
F_{i_1}(w_{\sig(1)}\ddd w_{\sig(i_1)})\odos
F_{i_k}(w_{\sig(n-i_k+1)}\ddd w_{\sig(n)}))=\\
\sum_{k=2}^n\sum_{I\in\NN^k\atop |I|=n}\sum_{\sig\in\Sigma_n}\frac{1}{I!k!}
\epsilon(\sig)
\sum_{l=1}^{k-1}\sum_{\tau\in\Sh(l,k)}
\epsilon(\sig,F_{i_1}(\Diamond)\ddd F_{i_k}(\Diamond))
F_{i_{\tau(1)}}(\Diamond)\odos
F_{i_{\tau(l)}}(\Diamond)\ot\\ 
F_{i_{\tau(l+1)}}(\Diamond)\odos
F_{i_\tau(k)}(\Diamond)=\\
\sum_{k=2}^n\sum_{I\in\NN^k\atop |I|=n}\sum_{l=1}^{k-1}{k \choose l}
\sum_{\sig\in\Sigma_n}\frac{1}{I!k!}
\epsilon(\sig,w_1\ddd w_n)
F_{i_1}(\Diamond)\odos
F_{i_l}(\Diamond)\ot\\ 
F_{i_{l+1}}(\Diamond)\odos
F_{i_k}(\Diamond)
\end{multline*}

Here, we have set $F_{i_m}(\Diamond):=
F_{i_m}(w_{\sig(i_1\kkk i_{m-1}+1)}\ddd w_{\sig(i_1\kkk i_m)})$.
Since ${k \choose l}\cdot\frac{1}{k!}=\frac{1}{l!(k-l)!}$, we see that
both sides coincide. 
Hence, by the commutativity of diagram (\ref{qudiag}), the difference
$$ F_n(w_1,...,w_n)-\sum_{k=2}^n\sum_{i_1\kkk i_k=n}\sum_{\sig\in\Sigma_n}
\frac{1}{I!k!}\epsilon(\sig)
F_{i_1}(w_{\sig(1)}\ddd w_{\sig(i_1)})\odos
F_{i_k}(w_{\sig(n-i_k+1)}\ddd w_{\sig(n)})$$
belongs to $\Kern(\Delta^+)=W'$. Thus it is just the term
$F_n(w_1,...,w_n)$, and the induction step is done.
\qed

\textbf{Proof of Proposition~\ref{Qcoal}:}
Induction on $n$.
Set $q_0:=Q_{0,1}(1)$.
(1) By the commutativity of diagram (\ref{qdelta}), comparing terms of
polynomial degree zero, we have that $Q_{i,0}=0$, for each $i\geq 0$.
(2) By the commutativity of diagram (\ref{qdelta}), comparing terms of
polynomial degree $i$ and linear degree $+1$, we have that
$Q_{0,i}=0$, for $i\neq 1$. (3) By similar arguments, we see that
for $w\in W$, we have $Q_{1,2}(w)=q_0\odot w$ and
$Q_{1,i}=0$, for $i\geq 3$. Thus the cases $n=0,1$ are done.\\
Now suppose that the statement is proved 
for all $m\leq n-1$.
Then $(Q\ot 1+1\ot Q)(\Delta(w_1\ddd w_n))$ can be written in the form
\begin{multline*}
\sum_{l=0}^{n-1}\sum_{k=1}^{n-l}\sum_{\sig\in\Sh(k+l-1,n)}
\sum_{\tau\in\Sh(l,k+l-1)}\epsilon(\sig)
\epsilon(\tau,w_{\sig(1)}\ddd w_{\sig(k+l-1)})
Q_l(w_{\sig(\tau(1))}\ddd w_{\sig(\tau(l))})\odot\\
\odot w_{\sig(\tau(l+1))}\odos w_{\sig(\tau(k+l-1))}\ot
w_{\sig(k+l)}\odos w_{\sig(n)}+\\
+\sum_{l=0}^{n-1}\sum_{k=1}^{n-l}\sum_{\sig\in\Sh(k,n)}
\sum_{\tau\in\Sh(l,n-k)}\epsilon(\sig)
\epsilon(\tau,w_{\sig(k+1)}\ddd w_{\sig(n)})
w_{\sig(1)}\odos w_{\sig(k)}\ot\\
\ot Q_l(w_{\sig(k+\tau(1))}\ddd w_{\sig(k+\tau(l))})\odot
w_{\sig(k+l+1)}\odos w_{\sig(n)}
\end{multline*} 
This is just the sum over k and l of the following expression:
\begin{multline*}
\sum_{\sig\in\Sh(l,k+l-1,n)}\epsilon(\sig)
Q_l(w_{\sig(1)}\ddd w_{\sig(l)})\odot
w_{\sig(l+1)}\odos w_{\sig(k+l-1)}\ot
w_{\sig(k+l)}\odos w_{\sig(n)}+\\
+\sum_{\sig\in\Sh(k,k+l,n)}\epsilon(\sig)(-1)^{w_{\sig(1)}\kkk w_{\sig(k)}}
\epsilon(\tau,w_{\sig(1)}\ddd w_{\sig(k+l-1)})
w_{\sig(1)}\odos w_{\sig(k)}\ot\\
\ot Q_l(w_{\sig(k+1)}\ddd w_{\sig(k+l)})\odot
w_{\sig(k+l+1)}\odos w_{\sig(n)}.
\end{multline*}
Here, by $\Sh(l,m,n)$ we mean the set of all permutations $\sig\in\Sigma_n$,
such that $\sig(1)<...<\sig(l)$ and $\sig(l+1)<...<\sig(m)$ and
$\sig(m+1)<...<\sig(n)$. 
On the other side,
\begin{equation*}
\Delta(\sum_{l=0}^{n-1}\sum_{\sig\in\Sh(l,n)}
\epsilon(\sig)
Q_l(w_{\sig(1)}\ddd w_{\sig(l)})\odot
w_{\sig(l+1)}\odos w_{\sig(n)}) 
\end{equation*}
can be written as sum over $k$ and $l$ of expressions of the form
\begin{eqnarray*}
\sum_{\sig\in\Sh(l,n)}\sum_{\tau\in\Sh(k,l-n+1)}\epsilon(\sig)
\epsilon(\tau,u_1\ddd u_{n-l+1})u_{\tau(1)}\odos u_{\tau(k)}\ot
u_{\tau(k+1)}\odos u_{\tau(n-l+1)},
\end{eqnarray*}
where we have set $$(u_1\ddd u_{n-l+1}):=(Q_l(w_{\sig(1)}\ddd w_{\sig(l)}),
w_{\sig(l+1)}\ddd w_{\sig(n)}).$$
We see easily that on both sides we have the same sums.
This finishes the induction step.
\qed

\textbf{Proof of Lemma~\ref{obst}:}
By our hypothesis, for each $m\leq n-1$, we have
\begin{equation}\label{hi}
Q_1\circ f_m-f_m\circ Q_1^{(m)}=
\sum_{k+l=m+1\atop k\geq 2}f_l\circ Q_k^{(n)}-\sum_{k=2}^m
\sum_{I\in\NN^k\atop |I|=m}Q_k\circ f_I.
\end{equation}
Furthermore, we can generalize the fact that $Q$ is an $L_\infty$-structure
to the the following equations: Let $m,k,l$ be natural numbers
such that $m\geq 1$ and $k+l=m+1$. Then
\begin{equation}\label{hii}
Q_1^{(l)}\circ Q_k^{(m)}+Q_k^{(m)}\circ Q_1^{(m)}+
\sum_{r+s=k+1\atop r,s\geq 2}Q_r^{(m+1-s)}\circ Q_s^{(m)}=0.
\end{equation}
Of course, they are also correct for $Q'$.
Now, we apply $\delta$ on the first summand of $r(f_1,...,f_{n-1})$.
Using equations (\ref{hii}) and (\ref{hi}), it takes the following form:
\begin{multline*}
\sum_{k=2}^n\sum_{I\in\NN^k\atop |I|=n}
(Q'_k\circ f_I\circ Q_1^{(n)} - Q'_1\circ Q'_k\circ f_I)=\\
\sum_k\sum_I Q'_k\circ{Q'_1}^{(k)}\circ f_I+
\sum_k\sum_I\sum_{r+s=k+1\atop r,s\geq 2}Q'_r\circ{Q'_s}^{(k)}\circ f_I+
\sum_K\sum_I Q'_k\circ f_I\circ Q_1^{(n)}=\\
-\sum_{k=2}^{n-1}\sum_I\sum_{\nu=1}^k\sum_{r=2}^{i_\nu}\sum_{j_1+...+j_r=i_\nu}
Q'_k\circ (1^{\ot\nu-1}\ot Q'_r\ot 1^{\ot k-\nu})\circ
f_{(i_1,...,i_{\nu-1},j_1,...,j_r,i_{\nu+1},...,i_k)}+\\
\sum_{\bar{k}=2}^n\sum_{|\bar{I}|\in\NN^{\bar{k}}\atop |\bar{I}|=n}
\sum_{\bar{r}+\bar{s}=\bar{k}+1\atop\bar{r},\bar{s}\geq 2}
\sum_{\bar{u}+\bar{v}=\bar{r}-1} Q'_{\bar{r}}\circ
(1^{\ot\bar{u}}\ot Q'_{\bar{s}}\ot 1^{\bar{v}})\circ f_{\bar{I}}+\\
\sum_{k=2}^n\sum_I\sum_{\nu=1}^k\sum_{s+t=i_\nu+1\atop t\geq 2}
Q'_k\circ f_{(i_1,...,i_{\nu-1},s,i_{\nu+1},...,i_k)}\circ
(1^{\ot i_1+...+i_{\nu-1}}\ot Q_t^{(i_\nu)}\ot 1^{\ot i_{\nu+1}+...+i_k}).
\end{multline*}
The first and second summand annihilate each other, since we have
the following 1:1 - correspondence of index sets:
$$(k,I,\nu,r,J)\mapsto
 (\bar{k}=k+r-1,\bar{I}=(i_1,...,i_{\nu-1},j_1,...,j_r,i_{\nu+1},...,i_k),
 \bar{r}=k,\bar{u}=\nu-1),$$
$$(k=\bar{r},I=(\bar{i}_1,...,\bar{i}_{\bar{u}},
  \bar{i}_{\bar{u}+\bar{s}+1},...,\bar{i}_{\bar{k}}), \nu+\bar{u}+1,
  r=\bar{s},J=(\bar{i}_{\bar{u}+1},...,\bar{i}_{\bar{u}+\bar{s}}))
\leftarrow (\bar{k},\bar{I},\bar{r},\bar{u}).$$

We apply $\delta$ on the second summand of $r(f_1,...,f_{n-1})$:
\begin{multline*}
\delta(\sum_{k=2}^n\sum_{I\in\NN^k\atop |I|=n}Q'_k\circ f_I)=
\sum_k\sum_I(Q'_k\circ f_I\circ Q_1^{(n)}-Q'_1\circ Q'_k\circ f_I)=\\
\sum_k\sum_I Q'_k\circ f_I\circ Q_1^{(n)}-
\sum_k\sum_IQ'_k\circ {Q'_1}^{(k)}\circ f_I-
\sum_k\sum_I\sum_{u+v=k+1\atop u,v\geq 2}Q'_u\circ{Q'_v}^{(u)}\circ f_I=\\
\sum_{k=2}^{n-1}\sum_I\sum_{\nu=1}^k\sum_{r=2}^{i_\nu}
\sum_{j_1+...+j_r=i_\nu} Q'_k\circ (1^{\ot\nu-1}\ot Q'_r\ot 1^{\ot k-\nu})
\circ f_{(i_1,...,i_{\nu-1},j_1,...,j_r,i_{\nu+1},...,i_k)}\\
-\sum_k\sum_I\sum_{s+t= i_\nu+1\atop t\geq 2}
Q'_k\circ f_{(i_1,...,i_{\nu-1},s,i_{\nu+1},...,i_k)}\circ
(1^{\ot i_1+...+i_{\nu-1}}\ot Q_t^{(i_\nu)}\ot 1^{\ot i_{\nu+1}+...+i_k})\\
-\sum_{\bar{k}}\sum_{\bar{I}}\sum_{u+v=\bar{k}+1\atop u,v\geq 2}
\sum_{c+d=u-1}Q'_u\circ (1^{\ot c}\ot Q'_v\ot 1^{\ot d})\circ f_I.
\end{multline*}

The first and third summand annihilate each other, since
we have the following 1:1 - correspondence of index sets:

$$(k,I,\nu,r,J)\mapsto (\bar{k}=k-1+r,
\bar{I}=(i_1,...,i_{\nu-1},j_1,...,j_r,i_{\nu+1},...,i_k), u=k,c=\nu-1)$$
and
$$(k=u,I=(\bar{i}_1,...,\bar{i}_c,\bar{i}_{c+v},...,\bar{i}_{\bar{k}}),
\nu=c+1,r=v,J=(\bar{i}_{c+1},...,\bar{i}_{c+v-1})).$$

The second term is just the remaining term above. So the statement
is proved.
\qed

\small{

\lz
\begin{center}
Institut Fourier\\
UMR 5582\\
BP 74\\
38402 Saint Martin d'H\`eres\\
France\\
\lz
frank.schuhmacher@ujf-grenoble.fr
\end{center}
}
\end{document}